\documentclass[default,iicol]{sn-jnl}

\jyear{2025}

\theoremstyle{thmstyleone}
\usepackage{url}
\usepackage{amssymb,amsthm,amsmath,amsfonts,amstext}
\usepackage{latexsym}
\usepackage{epsfig}
\usepackage[all]{xy}

\newtheorem*{thmm}{Theorem.}

\newtheorem*{proba}{Problems.}

\newtheorem*{proba'}{Problem A$\prime$.}

\begin{document}

\title[Unions of triangles]
{Regular triangle unions with maximal number of sides} 
\author{Giedrius Alkauskas}
\affil{\orgdiv{Institute of Computer Science, Department of Mathematics and Informatics}, \orgname{Vilnius University},  \orgaddress{\street{Naugarduko 24},
		\city{Vilnius}, \postcode{LT-03225}, \country{Lithuania}}. \\
	Corresponding author email: \texttt{giedrius.alkauskas@mif.vu.lt}}

\raggedbottom

\abstract{Fix an integer $n\geq 1$. Suppose that a simple polygon is the union of $n$ triangles whose vertices along the common boundary are arranged cyclically. How many sides can such a union -- to be called regular -- have at most? This gives OEIS sequence A375986, a recent entry. It will be shown here that the sequence begins $3, 12, 22, 33, 45, 56, 67, 80, 91$, and satisfies linear lower and upper bounds. The latter is not merely an estimate: it is realizable combinatorially. This leads to two further questions: can the same combinatorics be realized in pseudoline geometry, and if so, can such a realization be stretched? The paper is largely expository, with excursions into neighboring topics (union complexity, the Zone Theorem, stretchability, the Kobon triangle problem, Davenport-Schinzel sequences, lower envelopes of line segments). However, it adds a new tool tailored for studying regular unions; namely, triangulation shifts. In essence, this is a method to represent any such $n$-union by a triangulation of a regular $(n+1)$-gon and its dynamical mutation.}

\maketitle

\subsection*{Unions} 
A \emph{curvilinear polygon} is a simple polygon whose edges are smooth curves $f:[0,1]\mapsto\mathbb{R}^{2}$ (a hyperbolic polygon in the disc model provides a visual example). Consider the union of $n$ such polygons. Its complement will typically consist of finitely many holes together with one unbounded component. Each of these is a curvilinear polygon on its own, just with the last one viewed inside out. The total number of their edges is called the \emph{combinatorial complexity} (in dimension $1$) of the union. The asymptotic behavior of the complexity depends strongly on the maximal number of intersection points (say, $s$) that the boundaries of any two polygons in the union can have. \\
\indent Once $s>2$ intersection points are allowed, the combinatorics may become very intricate. Har-Peled \cite{peled} wrote a beginner-friendly introduction to the subject. A more advanced exposition by Agarwal, Pach and Sharir \cite{agarwal} contains a wealth of further ideas: directions towards applications (linear programming, robotics, molecular modelling, constructive solid geometry), a review of Davenport-Schinzel sequences (which will appear later in our paper), the study of unions of plane objects under additional restrictions (for example, keeping track of how often only $2$ intersection points occur among pairs of objects), unions of ``fat" or convex objects. Not to mention the abundance of inspiring pictures in that paper. For $s\geq 4$ the complexity becomes quadratic unless one restricts the geometry of the objects in question, or considers only the outer boundary of the union. \\
\indent Take $n$ (ordinary) triangles. Their sides will henceforth be called \emph{legs}, while \emph{edge} will be reserved for an edge of the union. Any two triangles intersect in at most $6$ points, and it is easy to contrive an example with $\sim 2n^2$ edges on the boundary of the union.  In such a construction, however, many of the triangles are extremely slim. If all triangles happen to have angles at least $\delta>0$ (so-called ``$\delta$-fat" triangles), then the complement of the union has at most $O_{\delta}(n)$ connected components, and the complexity is $O_{\delta}(n\log\log n)$ \cite{fat}.\\
\indent What, then, happens if the union has no holes at all and its interior is connected? Let the matter be put precisely. 

\begin{proba}
	Fix an integer $n\geq 1$. Suppose that a simple polygon is the union of $n$ triangles. What is the maximal number of sides it can have? The following variants arise naturally.
	\begin{itemize}
		\item[$\mathrm{1)}$]No further restrictions. This defines the sequence $i(n)$.
		\item[$\mathrm{2)}$]Vertices of triangles along the border of the union are arranged cyclically -- sequence $e(n)$.
		\item[$\mathrm{3)}$]Vertices of triangles lie on the unit circle -- sequence $c(n)$.
		\item[$\mathrm{4)}$]Vertices lie on the unit circle and are arranged cyclically -- sequence $r(n)$.
	\end{itemize}
\end{proba}
\noindent It may seem surprising that the question ``how many edges exactly?" does not seem to have been asked before. There are, however, reasons for this. The sequence $i(n)$ is the one to which existing methods and results apply most readily. Take a union of $n$ triangles with $i(n)$ edges. Find a triangle that contributes least to the common boundary. Remove it. It is not hard to see that this reduces the edge count of the outer border by at most $\frac{i(n)}{n}+3$. One is therefore tempted to write
\begin{eqnarray*}
	i(n-1)\geq i(n)-\frac{i(n)}{n}-3.
\end{eqnarray*}
If valid, this would imply $i(n)=O(n\log n)$. Of course, such a na\"{i}ve argument is not strictly correct: we do not know whether the removal of a single triangle leaves the union hole-free. Still, it captures part of the phenomenon. In the Section `Lower Envelopes', via the theorem of Micha Sharir \& Ady Wiernik (a deep and rather startling result to anyone new to the area), it will be seen that one cannot hope to do much better than this: $i(n)=\Omega(n\alpha(n))$, where $\alpha(n)$ is the functional inverse of the so-called (William) \emph{Ackermann function}. There are several variants of the latter. The one used here can be found in \cite{sharir}, p. 12--15. The exact variant is not important; all the asymptotic statements remain the same. The function $\alpha$ tends to $\infty$ with its argument, but at an extremely slow pace. Its value is $3$ throughout the interval $[5,16]$, and $4$ throughout the interval
\begin{eqnarray*}
	\Bigg{[}17, \underbrace{2^{2^{\ldots 2^{2}}}}_{65536\text{ copies}}\,\Bigg{]}.
\end{eqnarray*}
\begin{figure}
	\begin{center}
		\includegraphics[scale=0.20]{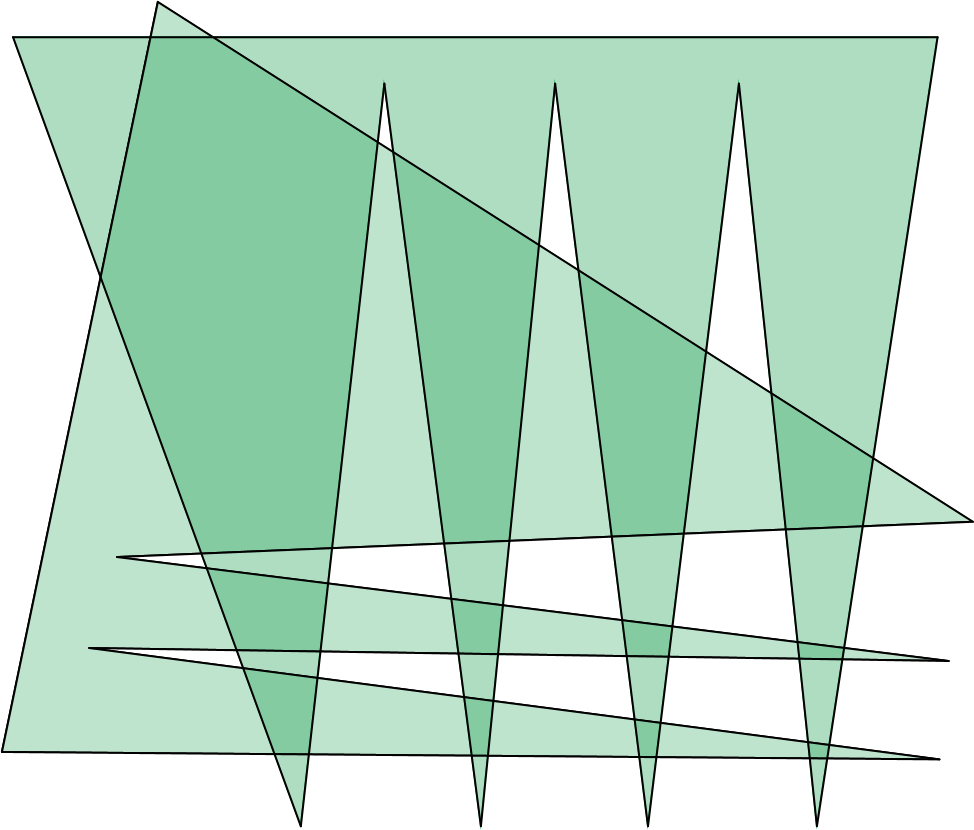} 
		\caption{$I=mn-m-n+3$ intersections of an $n$-gon and an $m$-gon for odd $m,n$ ($m=9$, $m=7$, $I=50$ here).} 
		\label{poly-two}
	\end{center}
\end{figure}
A sane mind is unlikely to picture where on the positive half-line this function first exceeds $5$. Hence no simple answer for $i(n)$ should be expected. The situation for $c(n)$ is even worse. Apart from the trivial inequality $c(n)\leq i(n)$, it is not clear where matters stand (see the last Section).\\
\indent The situation changes once the vertices are required to appear cyclically. Such unions will be called \emph{regular}. Unlike arbitrary unions, a regular one is automatically a simple polygon. Henceforth, $e(n)$ will be our main object of interest. It now appears in the OEIS as A375986 \cite{oeis}. The present paper, though largely expository, introduces regular unions for the first time and develops a tool particularly suited to them. In principle, this brings the exact value of $e(n)$ within reach for every $n$.\\
\indent There are, to be sure, problems in combinatorial geometry where even an arm's-length distance from a complete answer is still a substantial one. For example, it is known that a simple $m$-gon and an $n$-gon with an odd number of vertices can intersect in at least $mn-n-m+3$ points (Fig. \ref{poly-two}), and this is conjectured to be the largest possible. Eyal Ackerman et al. \cite{ackerman} proved that the maximum is indeed $mn-m-n+C$. Thus the problem is essentially settled, except for the exact value of $C$, which is still not known rigorously to equal $3$.\\
\indent What are the initial values of $e(n)$? For $n=1$, nothing needs proving. The union of $2$ triangles, if it is a simple $M$-gon, can have any number of sides in the range $3\leq M\leq 10$, and also $M=12$. Only in the latter case is the union regular. Metrically, the most familiar version is the regular hexagram (the Star of David). To prove that $e(3)=22$, a case-by-case analysis can be carried out. An example appears in Fig. \ref{startup} (three brown triangles).
\begin{figure}
	\begin{center}
		\includegraphics[scale=0.30]{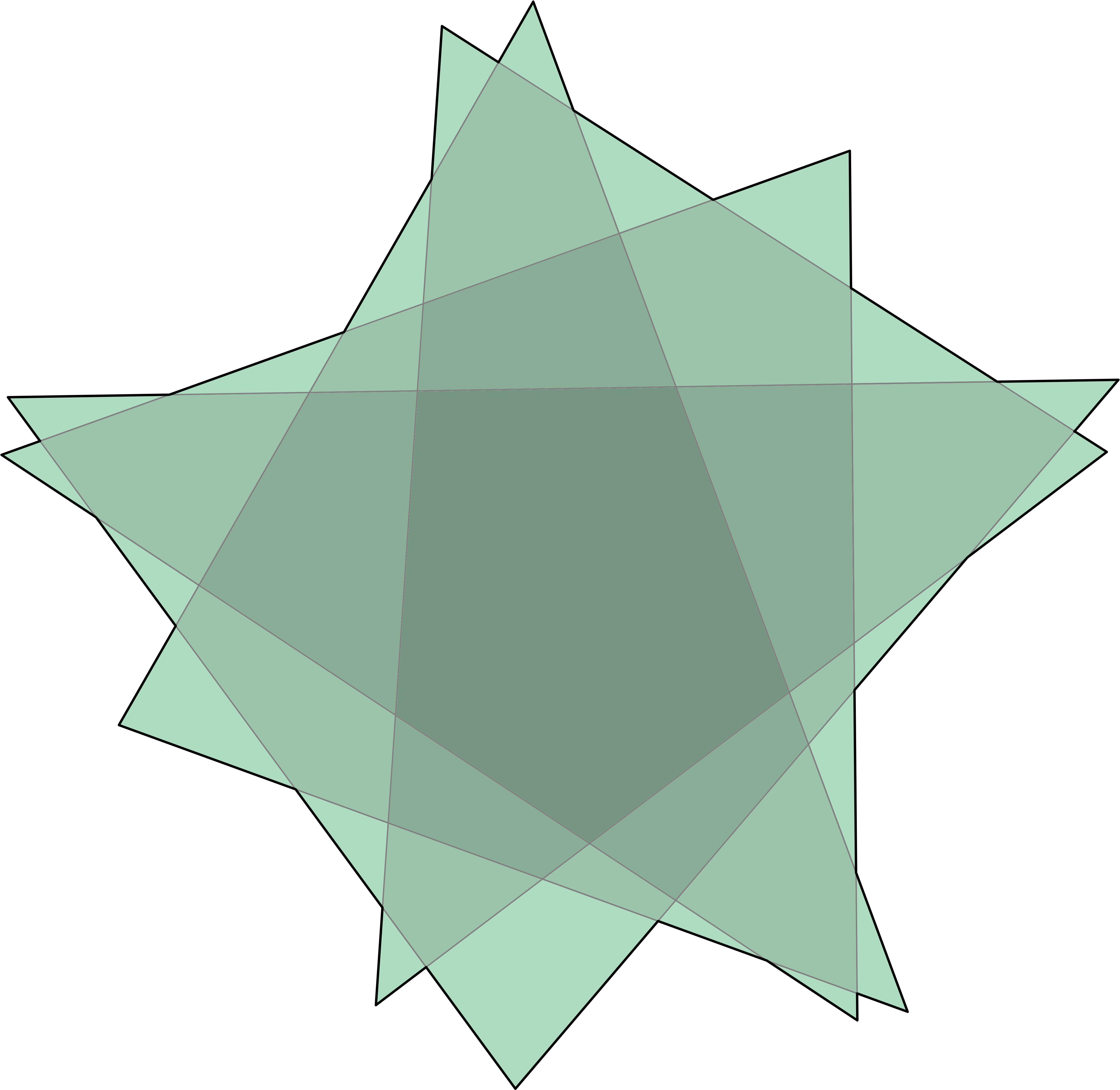} 
		\caption{Regular union of $4$ triangles with $33$ sides. This unique $33$-gon is calibrated to maximize the ratio of its area to the sum of squares of its edge-lengths.}
		\label{cont-4}
	\end{center}
\end{figure}
\begin{figure}
	\begin{center}
		\includegraphics[scale=0.32]{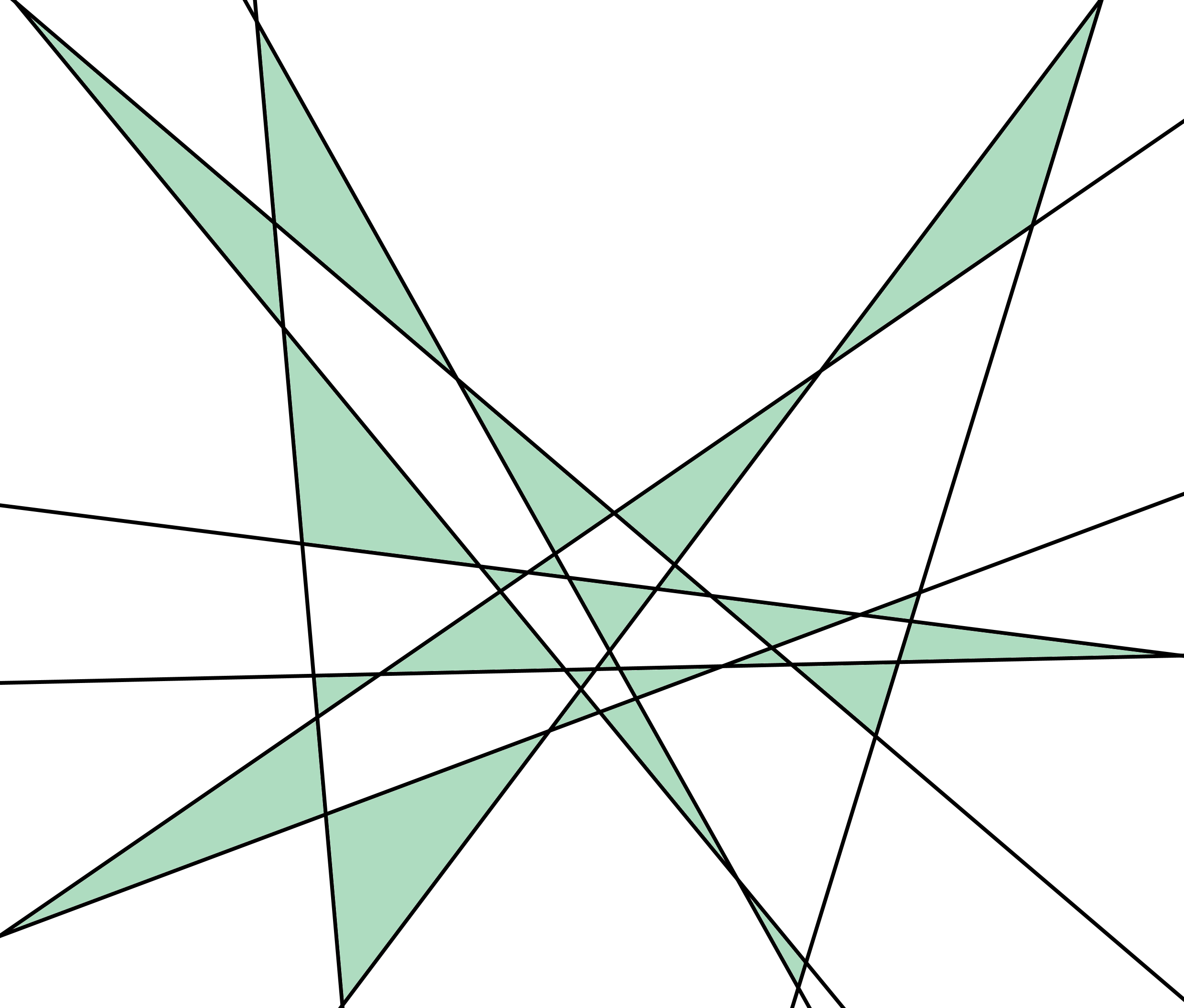} 
		\caption{$10$ lines forming $25$ non-overlapping triangles.} 
		\label{kobon}
	\end{center}
\end{figure}
\begin{figure*}
	\begin{center}
		\includegraphics[scale=0.63]{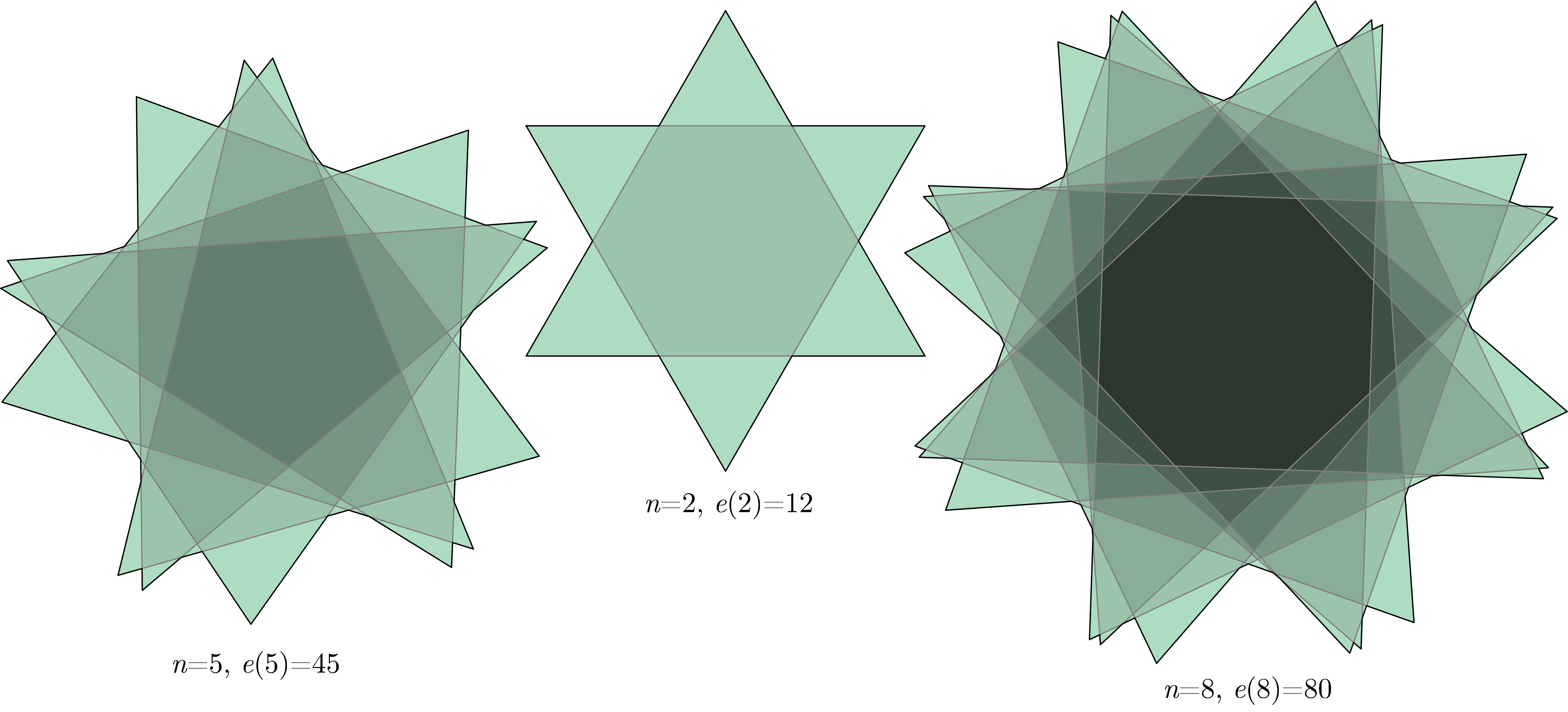} 
		\caption{All $3$ perfect proper unions of triangles. These (together with a single triangle) are exactly the unions satisfying $(n+1)M=(12n-6)n$. All are realizable on the unit circle.} 
		\label{tobuli}
	\end{center}
\end{figure*}
For $n=4$, a similar computer-assisted analysis shows that $e(4)=33$ (Fig. \ref{cont-4}). If all arrangements of $3n$ lines in general position in the affine plane (sequence A090338) were known, this catalogue would encode, among other things, the sequence $e(n)$ (and also $i(n)$). Such a list (for $n$ lines) would also encode the solution to the \emph{Kobon triangle problem} (sequence A006066), which asks how many non-overlapping triangles can be formed by $n$ lines (see Fig. \ref{kobon}, which depicts one of the maximal configurations for $n=10$). In full generality, Kobon's problem remains unsolved. In any case, a complete register of line configurations is of no practical use for large $n$, because their total number grows on the order of $2^{\Theta(n\log n)}$. The last exact value known is for $9$ lines, giving $3.111.341$ simple arrangements. In the present setting, that would settle, with some effort, the case $n=3$. Clearly, a different method is needed. In the setting of regular unions, such a method does emerge: \emph{triangulation shifts}. This is not some ingenious device. Rather, it appears naturally once the problem is examined in the right way. 
\begin{thmm} The sequence $e(n)$, $n\geq1$, starts as $3$, $12$, $22$, $33$, $45$, $56$, $67$, $80$, $91$. 
	For $n\geq 5$, the inequality
	\begin{eqnarray*}
		e(n)\leq 11n-14+2\bigg{\lfloor}\frac{n+1}{3}\bigg{\rfloor}=\mathscr{U}(n)
	\end{eqnarray*}
	holds. Additionally, $11\leq e(n+1)-e(n)\leq 14$.
\end{thmm}
\noindent The last result gives a situational lower bound $e(n)\geq 11n-8$ ($n\geq 10$). This can always be improved once we find a concrete $N\geq 10$ and a configuration of $N$ triangles with more than $11N-8$ edges. On the other hand, $\mathscr{U}(n)$ is not merely an estimate, but an exact combinatorial bound. The meaning of this phrase will become clearer as the argument unfolds. \\ 
\indent In the course of the proof, several excursions into neighboring topics will occur. The hope is that these do not obscure the picture. No complicated hidden scheme lies behind the proof, and the steps align transparently. As an intermediate result, it will follow that
\begin{eqnarray}
	e(n)\leq 12n-18+\frac{18}{n+1}.
	\label{pirm}
\end{eqnarray}	
This leads to the discovery of two ``perfect" unions and to a proof that $e(5)=45$, $e(8)=80$. Together with a triangle and a hexagram, this quartet forms the complete list of maximally saturated unions (Fig. \ref{tobuli}), where (\ref{pirm}) turns into equality. For $n\neq 1, 2, 5, 8$, the inequality must be strict.\\ 
\indent The value of $\mathscr{U}(n)$ stands in sharp contrast with the superlinear growth of the function $i(n)$. It is therefore worth pausing to ask what is meant by ``combinatorial" versus ``geometric" in describing this bound. That leads directly to pseudolines and stretchability. 

\subsection*{Pseudolines vs. Lines}
\begin{figure}[h]
	\begin{center}
		\includegraphics[scale=0.20]{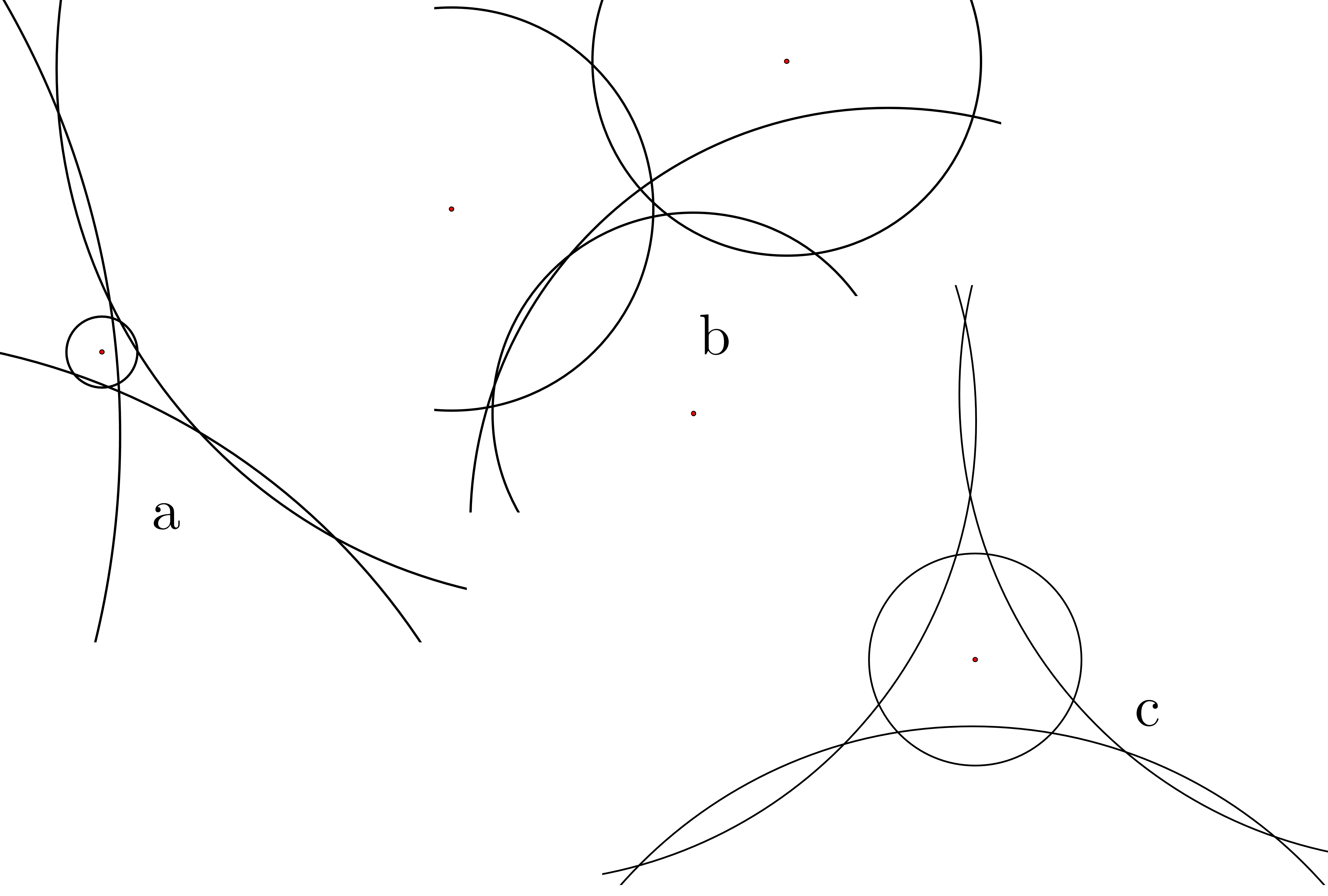} 
		\caption{Three ways to draw $4$ circles. Configuration $\mathbf{c}$ has $12$ boundary arcs.} 
		\label{dujen}
	\end{center}
\end{figure}

\begin{figure}
	\begin{center}
		\includegraphics[scale=0.35]{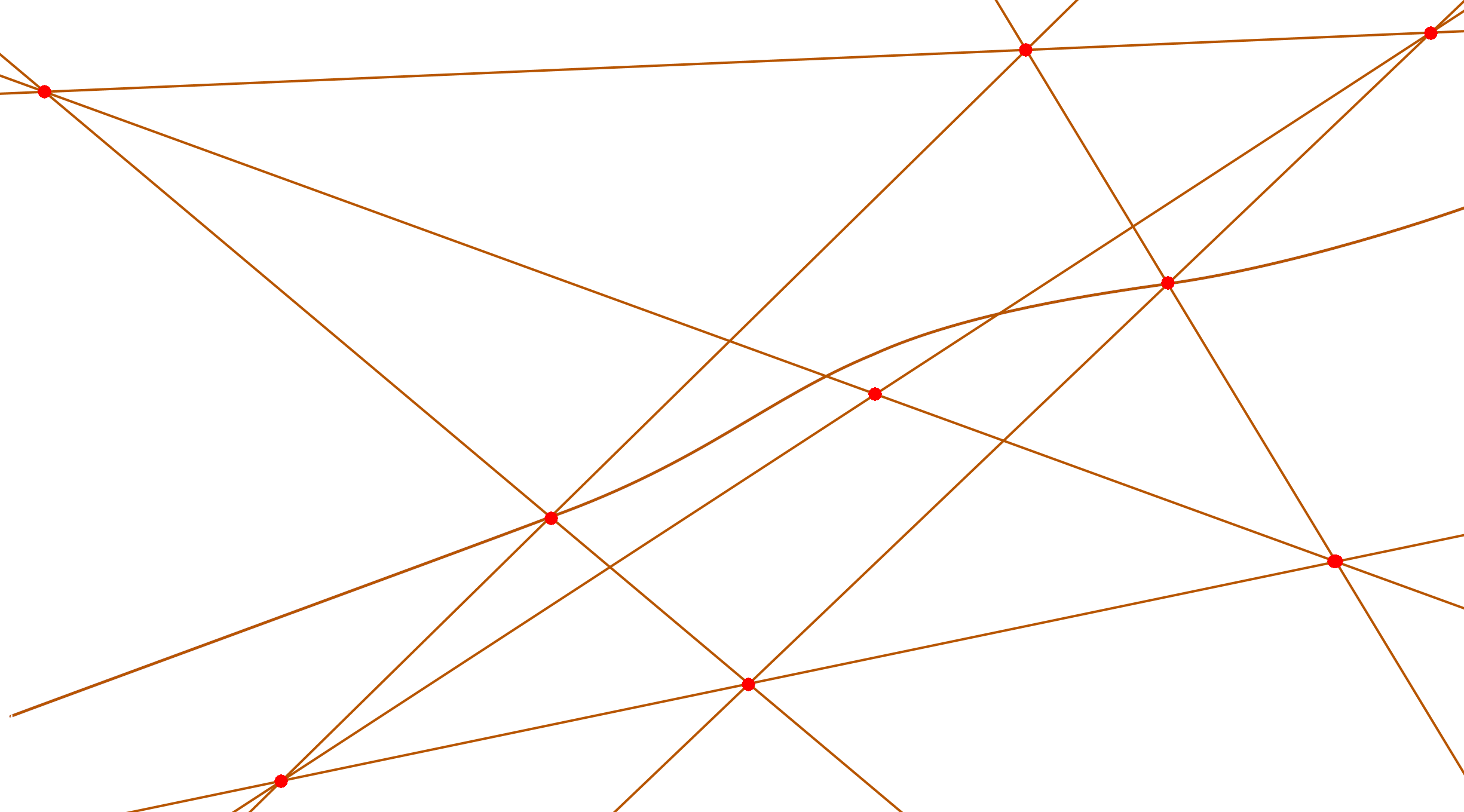} 
		\caption{Non-stretchable arrangement of $9$ projective pseudolines} 
		\label{pappus}
	\end{center}
\end{figure}

Besides A006066 and A090338, the OEIS contains other sequences made from kindred ingredients and carrying a similar flavor: sequences whose definitions use only familiar objects from Euclidean geometry, while their computation requires surprisingly delicate arrangements. One such sequence was popularized by Neil Sloane and the OEIS community during the ``ranking process" intended to decide which sequence should receive the honorary label A250.000. A runner-up in that vote, now sequence A250.001 (first proposed by Jonathan Wild), asks for the number of ways to draw $n$ circles in the affine plane \cite{sloane,oeis} (the notions of admissible and equivalent configurations are intuitive and can be found in those references). In Fig. \ref{dujen}, three of the $173$ arrangements of $4$ circles are reproduced.
See the recent survey \cite{cutler} for many further examples of such sequences, arising from cutting the plane with various ``exotic knives".\\
\indent Further still, problems such as pseudoline arrangements and stretchability are not only close in spirit and in delicacy, but they enter the present discussion in an essential way. A \emph{pseudoline} is a simple connected curve in the projective plane whose removal does not disconnect $\mathbb{R}P^{2}$. An \emph{arrangement} of pseudolines is a finite set of pseudolines such that every two meet exactly once, and the whole arrangement is not a pencil (there must be at least two intersection points). The affine (Euclidean) version is modified accordingly. A survey by Felsner and Goodman \cite{felsner} provides another rewarding point of entry. The most transparent arrangement of projective pseudolines that cannot be realized by straight lines -- that is, a non-stretchable one -- is the twisted Pappus configuration shown in Fig. \ref{pappus}. (The interpretation of Farey addition via Pappus's theorem even appeared on the cover of \emph{Intelligencer} \cite{schwartz}).
\begin{figure}
	\begin{center}
		\includegraphics[scale=0.18]{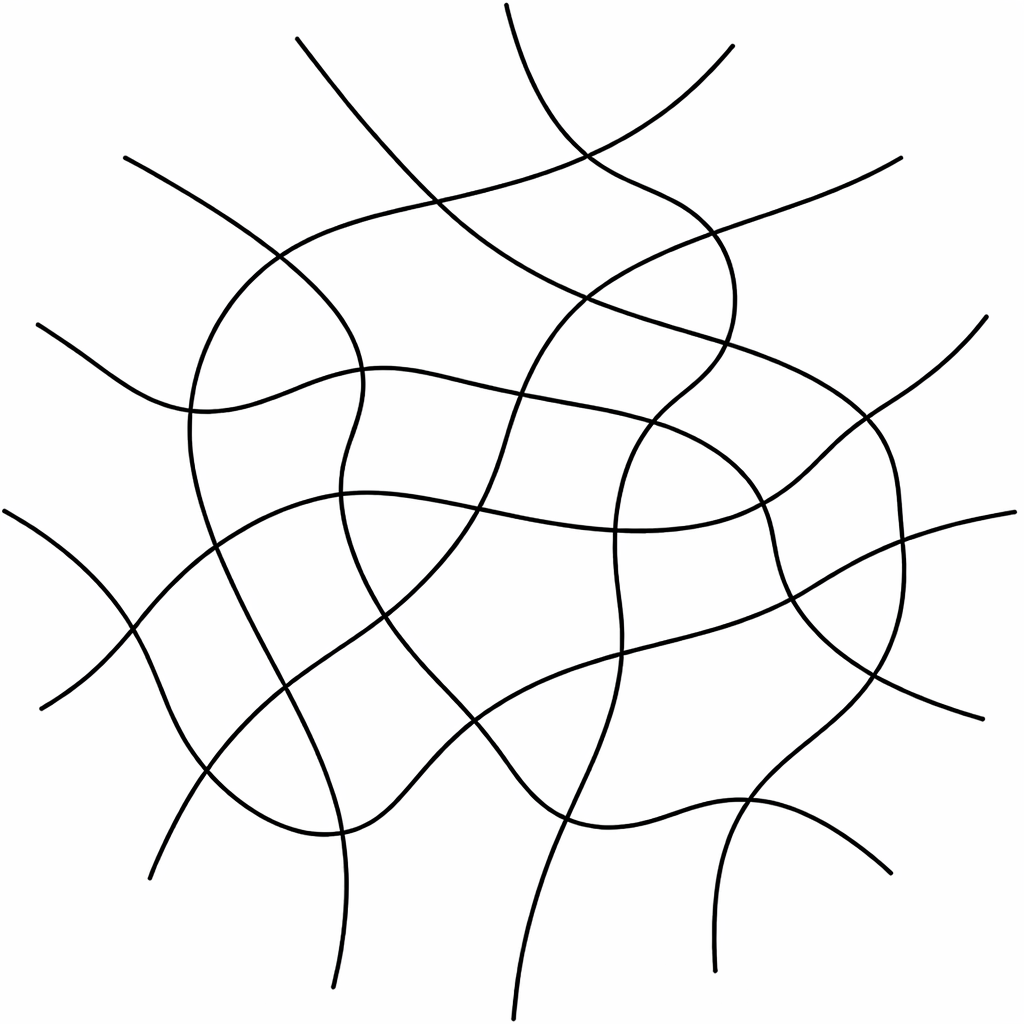} 
		\caption{One of $3$ non-stretchable simple arrangements of $8$ affine pseudolines; the rest ($38609$ in total) are stretchable.} 
		\label{pseudo}
	\end{center}
\end{figure}
The main object relevant to our paper is, however, slightly simpler: \emph{simple} arrangements of $n$ affine pseudolines (no multiple intersections allowed; sequence A090339).\footnote{There are some subtleties in what any given sequence actually counts. If $n$ affine lines are in general position, there are no finite multiple intersection points, but some subsets may be parallel. Thus, when viewed projectively, such an arrangement need not be simple. Other sequences count arrangements with a distinguished cell, with or without orientation, and so on.} For $n\leq 7$, every such arrangement can be realized by straight lines. Yet the configuration in Fig. \ref{pseudo} (found by L. Reeves and J. Wild) with $8$ pseudolines cannot be straightened. The best current bound on the number of simple arrangements of $n$ affine pseudolines, due to Justin Dallant, is $2^{0.6496n^2}$ \cite{dallant}. The proof relies on the \emph{Zone Theorem}, to which we return later. Compare this with the earlier bound for simple arrangements of $n$ straight lines; in light of this discrepancy, the aforementioned result of Wiernik and Sharir becomes even more striking (see the penultimate Section).\\
\indent Why introduce pseudolines at all? One reason is that their arrangement problems are easier, and some aspects can be handled in a purely combinatorial way. Many results about lines turn out to concern configuration more than straightness itself. There are further reasons. For example, Suvorov \cite{suvorov} gives two simple arrangements of $13$ lines in the affine plane. As pseudolines, one can be deformed continuously into the other while remaining simple throughout (\emph{flexible isotopy}). That is impossible if the curves are required to remain straight at every stage (\emph{rigid isotopy}). See the end of Chapter 4 in \cite{felsner} for a few more illuminating distinctions.\\
\indent Now take a regular union of $n$ triangles labelled $1$ through $n$, with $M$ edges. Walk once around the common boundary and record, in order, the label of the triangle to which each edge belongs. This way we obtain a cyclic sequence of length $M$. Much of what follows amounts to constructing a sequence of maximal possible length $\mathscr{U}(n)$ that serves as a model for such a boundary sequence. A second step would be to show that this model can indeed be realized in pseudoline geometry. That may well be manageable. Whether the resulting construction can then be straightened is a very different question. An analogous situation will appear later in the Section `Lower Envelopes'. The stretchability problem, as is known, is NP-hard. 

\subsection*{Lower Bound}
\begin{figure*}[h]
	\begin{center}	
		\includegraphics[scale=0.5]{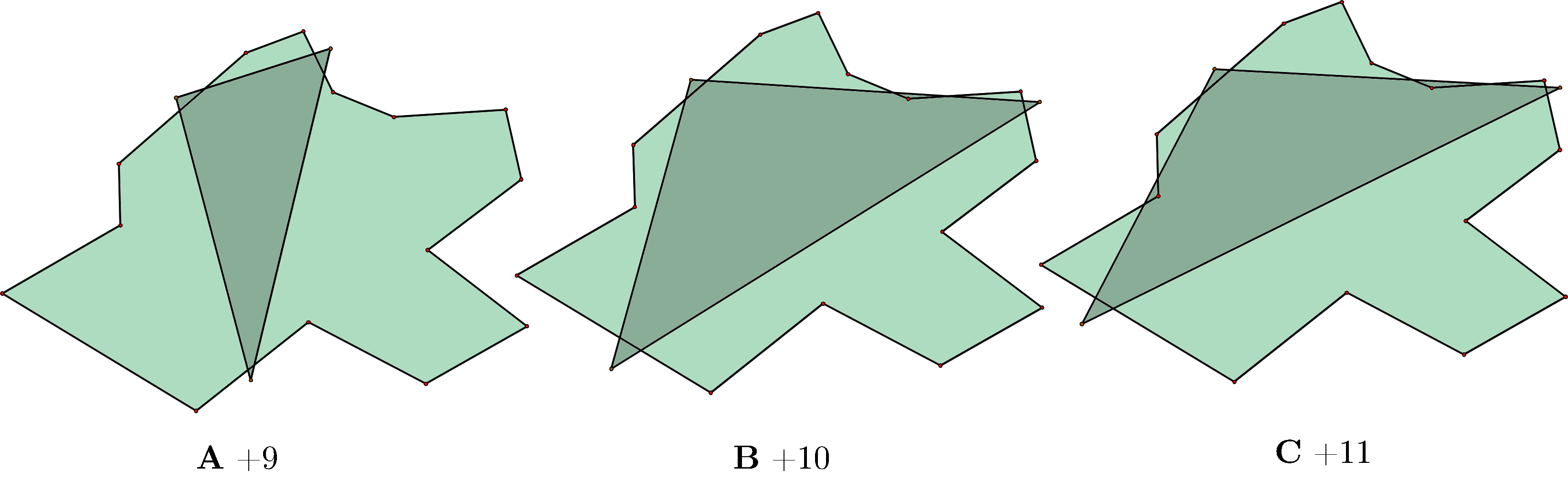} 
		\caption{Adding one triangle changes the edge count.}
		\label{by9}
	\end{center}
\end{figure*}

\begin{figure}[h]\begin{center}	
		\includegraphics[scale=0.27]{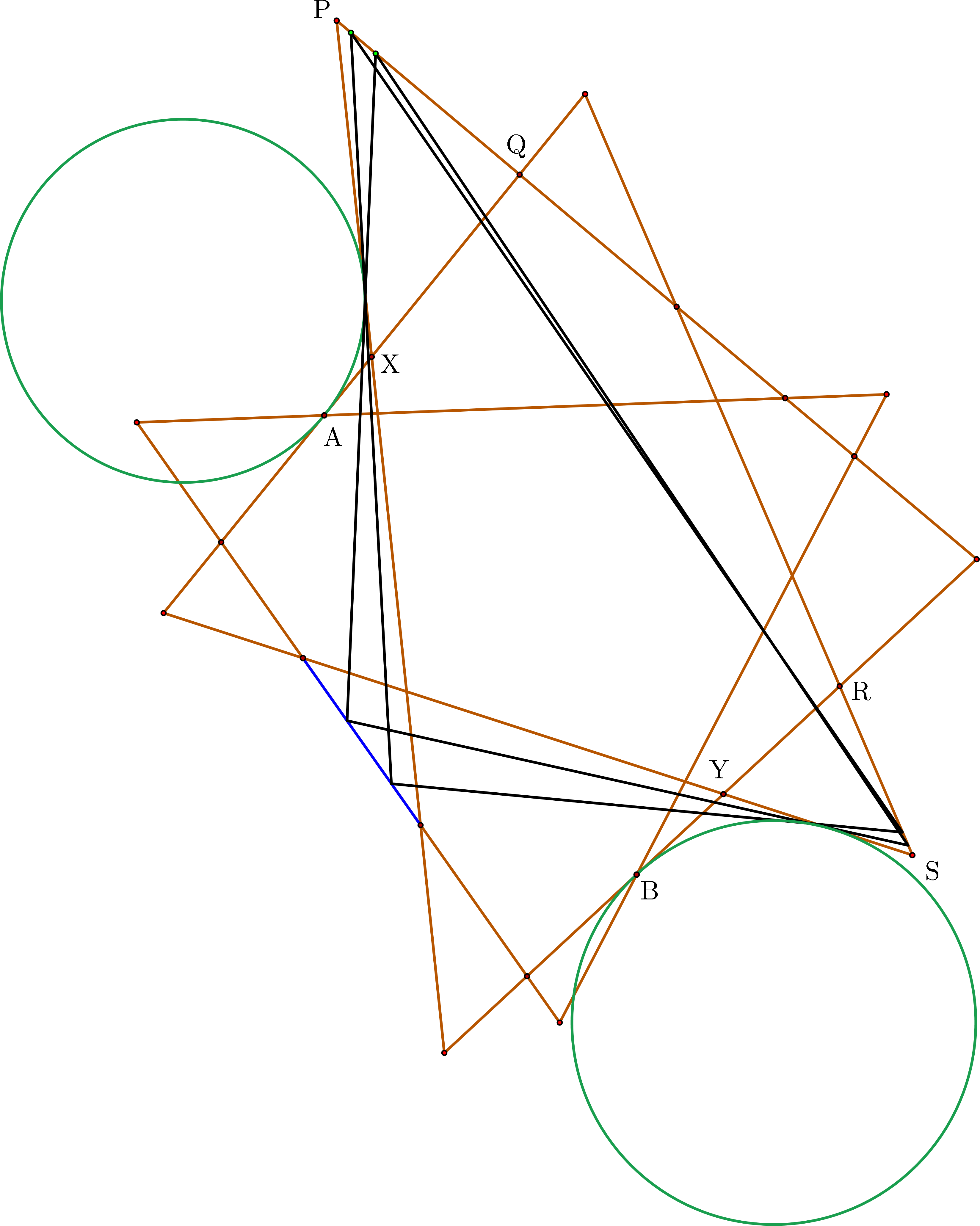} 
		\caption{An inductive construction showing that $e(n)\geq 11n-11$ for $n\geq 4$.}
		\label{startup}
	\end{center}
\end{figure}

\begin{figure*}[h]\begin{center}	
		\includegraphics[scale=0.64]{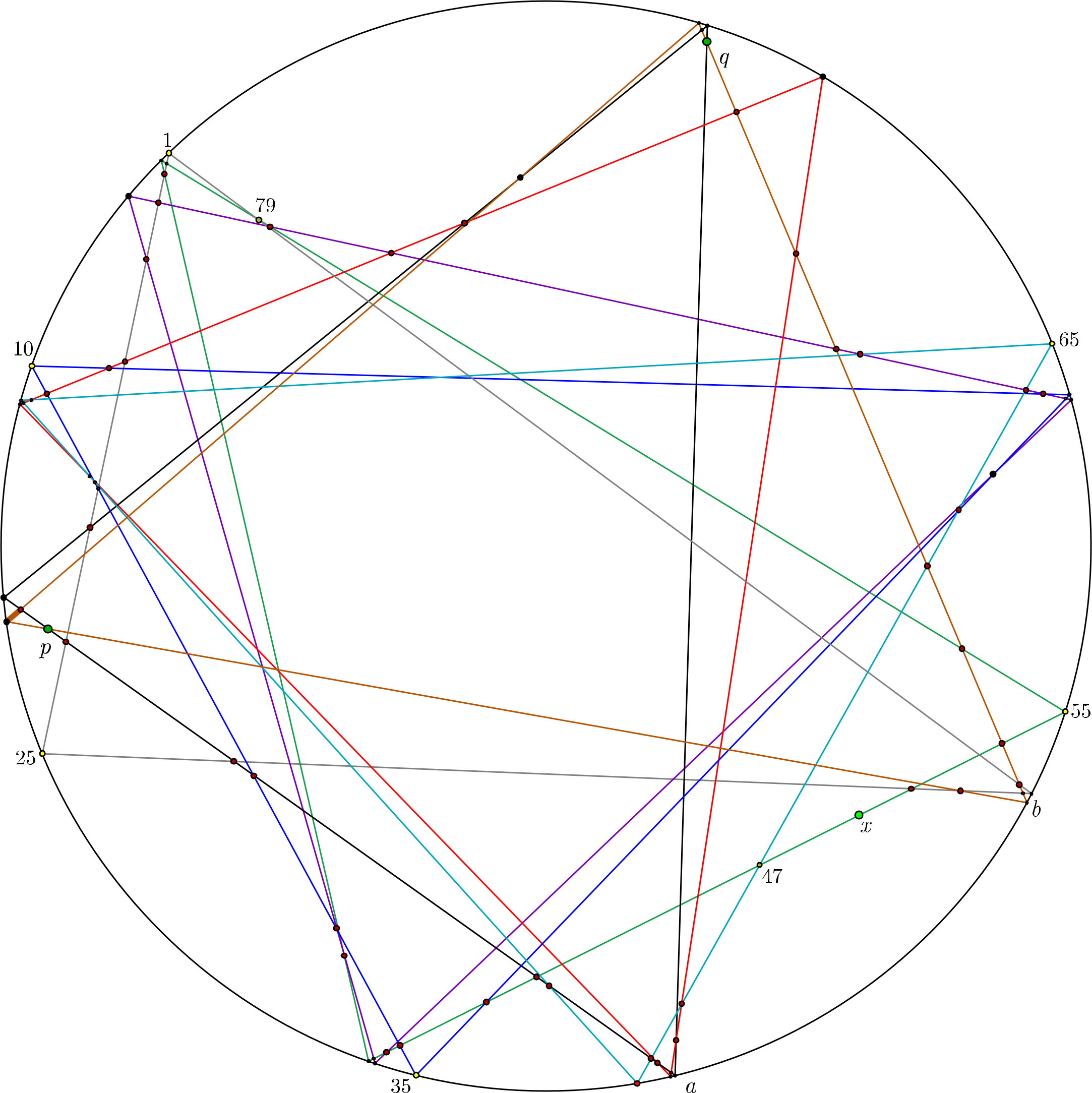} 
		\caption{A union of $8$ triangles with $79$ edges whose sequence is given by (\ref{sep-dev}). All vertices of the $8$ triangles lie on the unit circle; a few control vertices of the $79$-gon are labelled. Points $a,b,x,p,q$ are used in the inductive proof to show that $e(8+k)\geq 79+11k$.}
		\label{ses-sep}
	\end{center}
\end{figure*}

It is easy to construct, for any $n\geq 1$, a union of $n$ triangles that is a simple polygon with $9n-6$ edges. One simply chooses $3$ points on $3$ distinct sides of the existing polygon, pushes them slightly outward, and uses them to build a new triangle (Fig. \ref{by9}\textbf{A}).
The only two things to watch are that no vertex of the original polygon should lie inside the new triangle, and that the new polygon should remain simple. If, however, two incident edges are cut by the same leg of the new triangle, the edge count increases by $10$ (Fig. \ref{by9}\textbf{B}). If two disjoint pairs of incident edges are cut, the increase is $11$ (Fig. \ref{by9}\textbf{C}).\\ 
\indent With this in mind, let us start from the three brown triangles shown in Fig. \ref{startup}. Their union has $22$ edges. Consider the blue segment and the two green circles indicated there. Next choose several points, labelled $1,2,\ldots,m=n-3$ from left to right, strictly inside the blue segment. From each such point $p$, draw the inner tangents to both circles. Let their intersections with $PQ$ and $SR$ be $p'$ and $p''$, respectively. Fig. \ref{startup} shows the case $m=2$. 
Now push all points on $PQ$, $SR$, and on the blue segment slightly outward. Each of the $m$ new triangles contributes $11$ edges. This gives the total count $11n-11$.\\
\indent In fact, there is nothing special about this starting configuration. The same construction works whenever the union has an edge that is not adjacent to any original triangle vertex (hence $n\geq 3$). For example, take the green edge in Fig. \ref{ses-sep}, the one with control vertex ``47". Its endpoints belong to the cyan and grey triangles, respectively. Choose a point $x$ in the interior of that edge. A new triangle will be built with one vertex at $x$, pushed outward by just a bit. It will lie between the black and brown vertices $a$ and $b$. Let the left and right boundary edges at those vertices meet at points $p$ and $q$ (also indicated in Fig. \ref{ses-sep}). Draw the lines $xp$ and $xq$, and let them meet the corresponding brown and black edges of the union at $y$ and $z$. The slightly inflated triangle $xyz$ contributes $+11$ to the edge count. Its vertices also preserve the cyclic order. In this way the same mechanism as in Fig. \ref{startup} is reproduced.
\subsection*{Upper Bound}
Take a regular union of $n$ triangles with $M$ edges. For a concrete example, consider the $8$ triangles in Fig. \ref{ses-sep}. Start with the blue triangle (labelled $1$) and with the edge just before the control vertex ``10". The produced edge sequence is of length $79$:
\begin{eqnarray}
	\underline{\underline{11}\textcolor{red}{3}\overline{2233\textcolor{red}{216}4455\textcolor{red}{4}66\textcolor{red}{418}7788
			\textcolor{red}{7}11}}\overline{\textcolor{red}{74}22}\textcolor{red}{4}3344\textcolor{red}{3276}\nonumber\\
	5566\textcolor{red}{5}77\textcolor{red}{521}88\underline{\underline{11\textcolor{red}{8}22\textcolor{red}{85}33
			\textcolor{red}{5} 4455\textcolor{red}{4387}6677\textcolor{red}{6}88\textcolor{red}{63}}} \label{sep-dev}
\end{eqnarray}
\noindent Now isolate a fragment of this boundary lying between two vertices of the same blue triangle -- for example, the singly underlined one in (\ref{sep-dev}). Write down the numbers of the triangles whose RIGHT edges are encountered along the way. In other words: ignore a symbol until its double occurs; then record it once; write down every subsequent copy. In this way we obtain $\mathcal{L}=123214546417871$. What is intrinsic to this sequence? First, no two adjacent symbols are equal. Second, although symbols may appear more than once, no subsequence of the form $abab$ can occur for any pair $a,b\in\{1,2,\ldots,n\}$. This is a crucial feature of a regular triangle union. In the section `Lower Envelopes' it will be seen that once subsequences of the form $abab$ are allowed, the length may become superlinear.\\
\indent Why is this pattern missing? Since the boundary is traversed counterclockwise, each leg has a direction. After extending a leg to a full line, one may speak of its right and left sides. The same applies to any fixed edge. Every other edge then lies either to the right, to the left, or collinear with it. For any given point on the leg, one may also distinguish what lies behind and what lies ahead. Now fix one edge -- say one of the numbers in $\mathcal{L}$ -- and look at the subsequent edges in the right-edge sequence. They are either collinear with it or lie to the right of it. Suppose a subsequence $aba$ were to occur. Let the legs labelled $a$ and $b$ meet at a point $P$. As the boundary is traced, one passes from $a$ to $b$ and then back to $a$. At that moment the point $P$ lies behind. Returning to $b$ is impossible, because all points of $b$ ahead of $P$ lie to the left of $a$. \\
\indent Thus $\mathcal{L}$, and in fact every sequence of right edges between two same-colored vertices of a fixed triangle, is by definition a $(2,n)$-\emph{Davenport-Schinzel} sequence \cite{sharir}. The maximal length of such a sequence is known to be $2n-1$. The proof is elementary and inductive. Assume it for $n-1$, and consider a $(2,n)$-DS sequence $\mathcal{J}$. If not all symbols $1,\ldots,n$ occur, then $\mathcal{J}$ is a $(2,m)$-DS sequence for some $m<n$. Otherwise, reading from left to right, suppose the symbols first appear in the order $1,2,\ldots,n$ (renaming is always possible). Then the symbol $n$ can occur only once. Indeed, if it occurred twice, then some subsequence $npn$ would appear, and because of the chosen first-occurrence order, a subsequence $pnpn$ would be seen -- a contradiction. Remove $n$, and if its two neighbours are equal, merge them. What remains is a $(2,n-1)$-DS sequence of maximal length $2n-3$. Hence $\mathcal{J}$ has length at most $2n-3+2$. In a maximal $(2,n)$-DS sequence, incidentally, the last symbol to appear must always be flanked by two equal symbols.\\
\indent Return now to the edge count. Take the same underlined fragment, but read it from the opposite end. This time write down the triangles whose LEFT edges contribute to the boundary. One obtains the sequence $18781654613231$, again a $(2,n)$-DS sequence, though not a maximal one. It follows that the total number of edges between two blue vertices is at most $4n-2$, with the convention that the blue ones have been counted twice. Perform the same count for each of the three pairs of blue vertices. If $\ell_i$ denotes the number of edges contributed by the $i$-th triangle to the boundary of the $M$-gon, then
\begin{eqnarray*}
	2\ell_{1}+\ell_{2}+\cdots +\ell_{n}\leq 12n-6.
\end{eqnarray*}
Summing over all $i$ yields $(n+1)M\leq n(12n-6)$, which is precisely (\ref{pirm}). For the known values $n\leq 9$, this is already satisfactory, except that for $n=7$ the bound gives $68$, whereas only $67$ has so far been achieved $(45+11+11)$.\\
\indent The inequality $e(n)<12n$ has one immediate consequence. Take a regular union of $n$ triangles with $e(n)$ edges. Choose a triangle that contributes least to the common boundary. Its toll is then at most $11$ edges. Remove it. Each edge not adjacent to one of its vertices decreases the edge count by at most $1$, while two edges meeting at any such vertex may together decrease it by at most $3$. Hence $e(n-1)\geq e(n)-14$. \\
\indent The idea of separating left and right edges between two same-colored vertices was inspired by the proof of the Zone Theorem due to Edelsbrunner et al. \cite{edelsbrunner}. In that paper, however, `right' and `left' refer to two directions along the same ray.
\begin{figure*}[h]\begin{center}
		\includegraphics[scale=0.45]{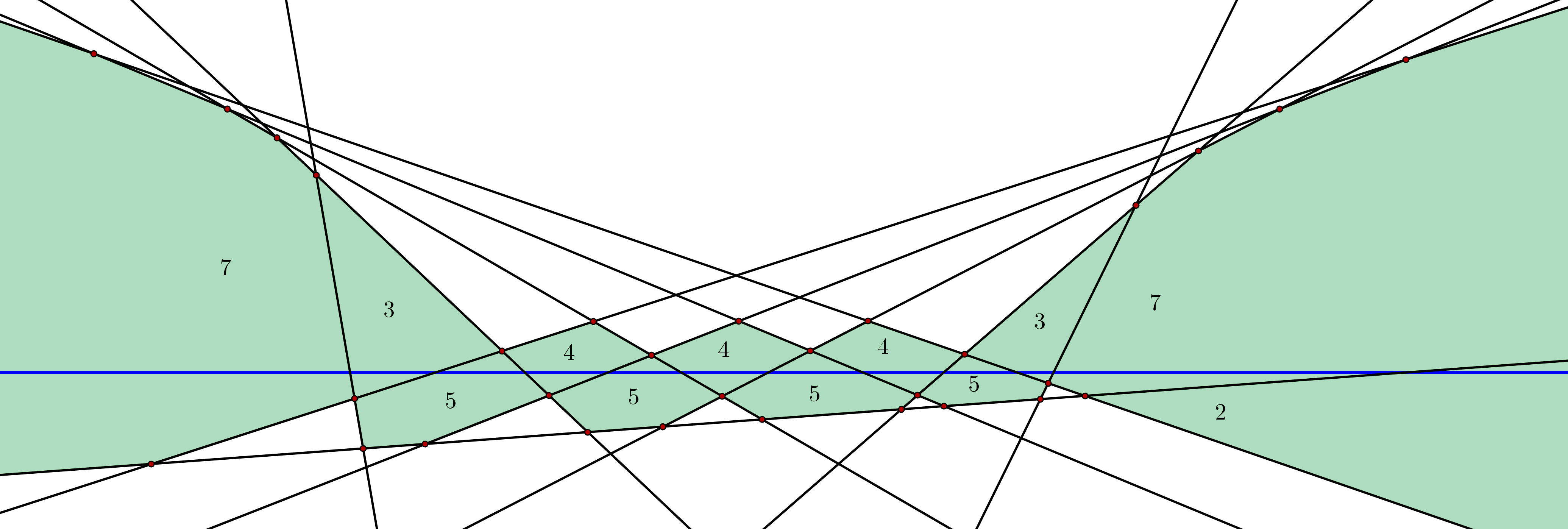} 
		\caption{An arrangement of $11$ lines; the complexity of the zone (colored in green) of a blue line is $54=\lfloor 5.5\cdot 11\rfloor-6$; the size of each face in the zone is marked.}
		\label{vienuol}
	\end{center}
\end{figure*}
Among the several proofs of the Zone Theorem, this one is notable for the essential use it makes of $(2,n)$-DS sequences.\\
\indent What does the theorem say? Let $\mathcal{K}$ be a simple arrangement of $n$ lines in the affine plane. The set $\mathbb{R}^{2}\setminus\mathcal{K}$ splits into connected components called \emph{faces}, each either a polygon (if bounded) or a region with two infinite edges. The \emph{size} of a face is the number of its edges. Given a line $\ell$, if $\mathcal{K}\cup\{\ell\}$ remains a simple arrangement, then the \emph{zone} of $\ell$ is the collection of all faces that it crosses. The \emph{complexity} of the zone is the sum of the sizes of those faces. This notion was introduced by Chazelle et al. in \cite{zone}, where an upper bound of $6n$ was established. The DS-sequence argument gives $6n-4$ (see \cite{edelsbrunner}, p. 334, and also \cite{sharir}, p. 125). In \cite{horizon}, Bern et al. improved this to $\lfloor 5.5n\rfloor-3$, and later Pinchasi \cite{pinchasi} sharpened it further to $\lfloor 5.5n\rfloor-5$. The arrangement in Fig. \ref{vienuol} (taken from \cite{horizon}, with a slight enhancement) shows that this last bound is nearly optimal. The similarity with Fig. \ref{startup} -- and later with Fig. \ref{two-cyc} -- is hard to miss. Such fan constructions, surely, are a basic tool employed throughout all the research that deals with unions.  

\subsection*{Maximal Sequences}
\begin{figure}
	\begin{center}
		\includegraphics[scale=0.28]{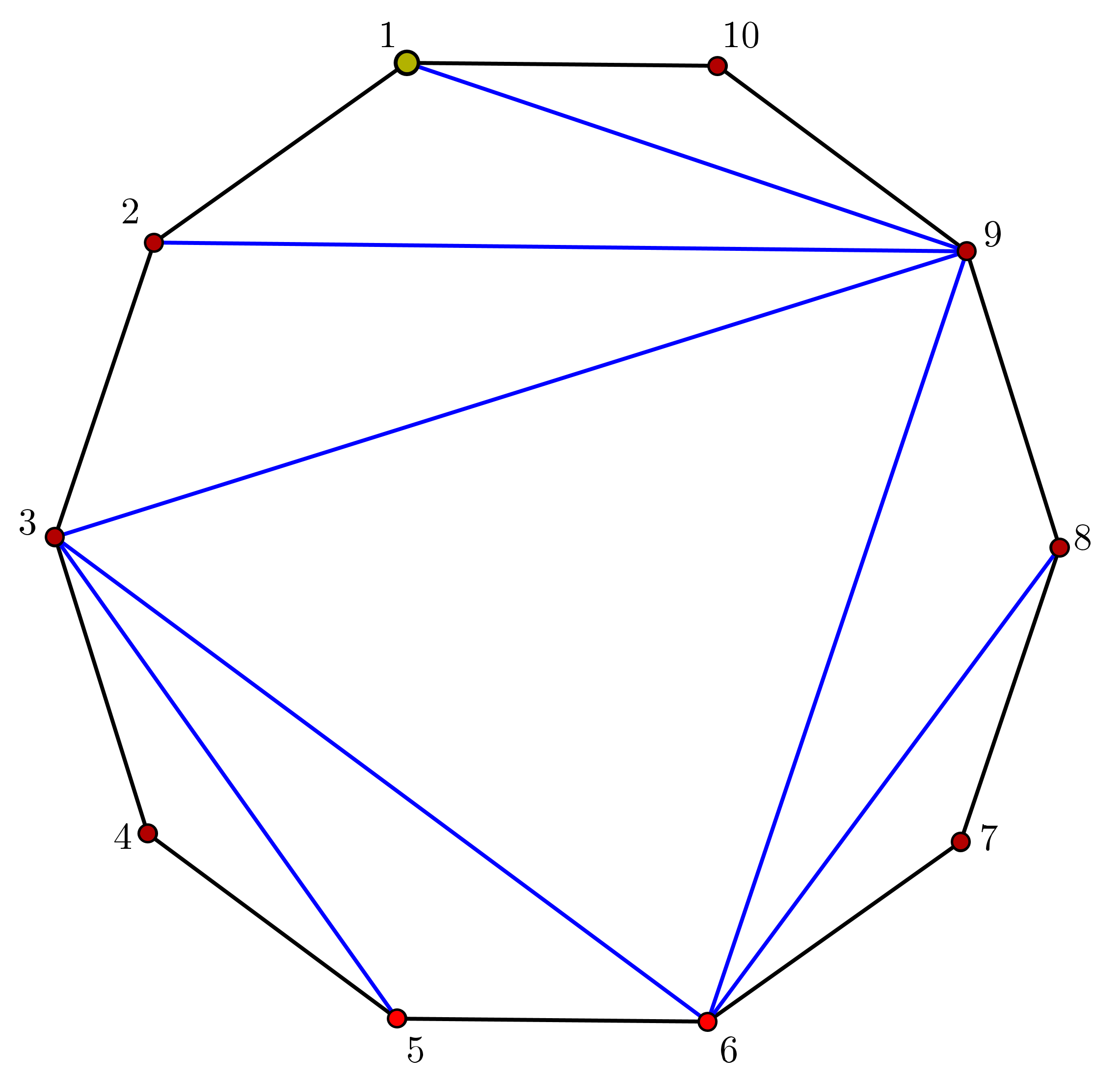} 
		\caption{A triangulation of a regular $10$-gon corresponding to the maximal $(2,9)$-DS sequence $12343536768632191$.}
		\label{pat}
	\end{center}
\end{figure}

The last refinement of the upper bound comes from observing that, depending on $n$, a certain number of the resulting $6n$ sequences must fail to be maximal. Once the related analysis is performed, one also gets an algorithm for constructing a sequence of length $\mathscr{U}(n)$ that respects the combinatorics of a regular union of $n$ triangles whose sides are segments of pseudolines. \\
\indent Begin again with $\mathcal{L}=123214546417871$. This is a maximal sequence. Now shift attention to the cyan triangle -- in other words, to the overlined fragment in (\ref{sep-dev}). The sequence extracted there is also maximal. Surprisingly, it can be reconstructed from $\mathcal{L}$ itself. Indeed, remove all $1$'s from $\mathcal{L}$ except the last one, and append a $2$ at the end. This gives $232454647871\ldots 2$.
For this to be maximal, the dots must conceal exactly two symbols, so the new sequence has the form $\mathcal{L}_{2}=232454647871ab2$. The claim is that $a$ and $b$ are forced. The reasoning will now certainly look familiar: let us methodically remove the last (reading from left to right) symbol to make its appearance, and merge its neighbours. Such a merge must occur:
\begin{eqnarray*}
	2324546478\underline{71a}b2&&\mathop{\longrightarrow}^{a=7}\,\,
	23245464\underline{787}b2\rightarrow\\
	2324546\underline{47b}2&&\mathop{\longrightarrow}^{b=4}\,\,232454642.
\end{eqnarray*}
\noindent Thus the only possible reconstruction is $\mathscr{S}(\mathcal{L})=\mathcal{L}_{2}=232454647871742$. The overlined fragment in (\ref{sep-dev}) confirms exactly this. Iterating the same procedure (the reconstructed symbols are printed in red) gives
\begin{eqnarray}
	&&\scriptstyle
	123214546417871\rightarrow232454647871\textcolor{red}{74}2\rightarrow3454647871742\textcolor{red}{4}3\rightarrow\nonumber\\
	&&\scriptstyle
	454647871742434\rightarrow5678717234\textcolor{red}{3276}5\rightarrow678717234327656\rightarrow\nonumber\\
	&&\scriptstyle
	7871723432756\textcolor{red}{5}7\rightarrow 81234325657\textcolor{red}{521}8\rightarrow\nonumber\\
	&&\scriptstyle
	123432565752181\rightarrow 2343256575281\textcolor{red}{8}2.\label{qr}
\end{eqnarray}
This agrees exactly (!) with the actual sequence extracted from (\ref{sep-dev}). So far, every one of them is maximal. Another pattern also appears: $\mathscr{S}^{9}(\mathcal{L})$ is obtained from $\mathcal{L}$ simply by adding $1$ to every term. This is a general rule.\\
\indent Let us stop iterating abstractly and return to the actual sequences drawn from (\ref{sep-dev}). For a while, one maximal sequence follows another. But at step $17$ (the doubly underlined fragment), a length $14$ sequence appears: $12345436768\textcolor{red}{63}1$. It is then followed by the maximal $2345436768631\textcolor{red}{3}2$.
The next two sequences are again maximal. Then a second shorter sequence turns up, $56768612321645$ (length $14$), followed by the maximal $6768612321645\textcolor{red}{4}6$. After two further steps everything closes up at the original $1232145464178\textcolor{red}{7}1$. Thus the process is almost smooth: over these $24$ steps only two glitches occur, and otherwise the next sequence is determined. The total number of edges in the union is the number of red symbols added, plus $6\cdot 8$, the number of black ones.\\
\indent If this already feels too opaque, that impression is valid. The strings in (\ref{qr}) look vague enough. Yet employing the help of a visualization of maximal $(2,n)$-sequences in terms of triangulations of a regular $(n+1)$-gon (first discovered by D. Roselle \cite{roselle}; see also \cite{sharir}, p. 250--232), something pleasant happens.\\
\indent To see this, consider the maximal $(2,9)$-DS sequence $12343536768632191$.
Split it into $9$ decreasing blocks: $\vert\,1\vert\, 2\,\vert\,3\,\vert\,43\,\vert\,53 \,\vert\,6\,\vert\,76 \,\vert\,86321\,\vert\,91\,\vert$. Now take a regular $(n+1)$-gon. For each vertex $i+1$, join it to all vertices that appear in the block beginning with $i$. This produces Fig. \ref{pat}. The picture is not merely suggestive: it gives a bijection between triangulations of an $(n+1)$-gon with a marked vertex and essentially distinct maximal $(2,n)$-DS sequences (that is, sequences not made equivalent by renaming the symbols).
\begin{figure*}
	\begin{center}
		\includegraphics[scale=0.82]{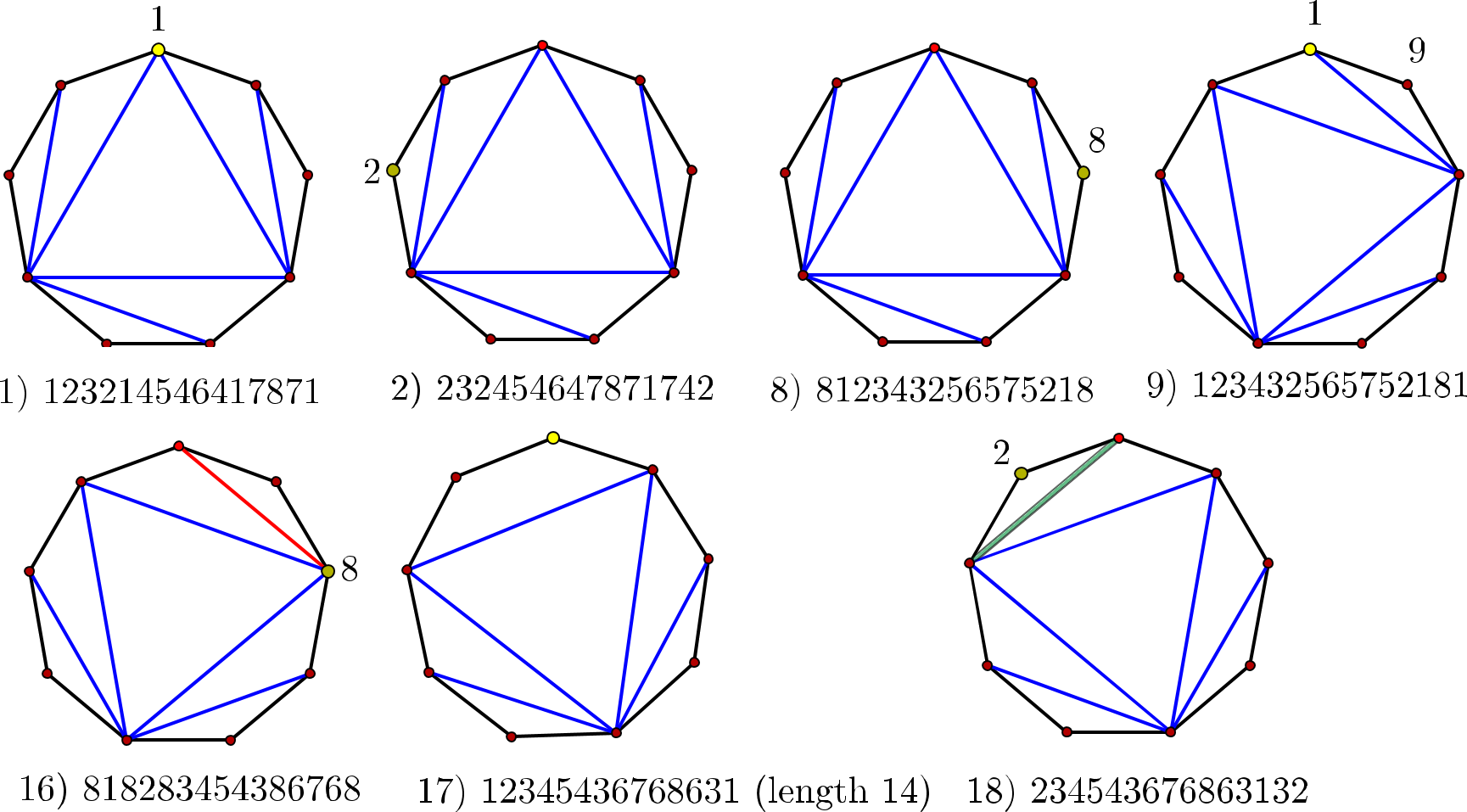} 
		\caption{Triangulations corresponding to the actual sequences of RIGHT edges in Fig. \ref{ses-sep}; 1--9 also correspond exactly to the sequences in (\ref{qr}).}
		\label{ratt}
	\end{center}
\end{figure*}
\begin{figure*}
	\begin{center}
		\includegraphics[scale=0.35]{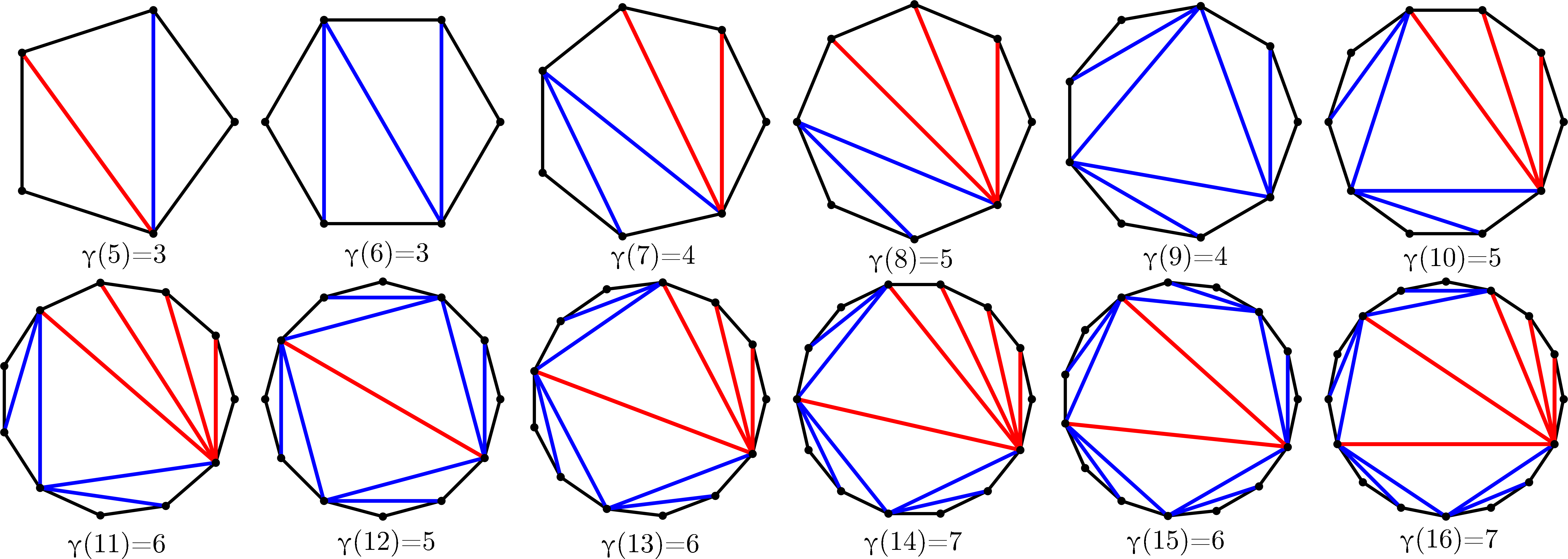} 
		\caption{Initial triangulations $\mathcal{T}_{0}$ of shifts with maximal weight $\omega(\mathcal{T}_{0})=12(n-1)-\gamma(n+1)$ (though weight is vertex-dependent, all vertices, save for two in cases $3\nmid n$, can be taken as initial vertex $1$). Red diagonals are those not present in the same triangulation, if rotated by $-3$ turns.} 
		\label{gamma}
	\end{center}
\end{figure*}
In the present setting this visualization convention is used with the natural order $1,2,\ldots,n$. If instead the order is $i,i+1,\ldots,i-1$, then the rule is as follows. First rotate the symbols cyclically until the natural order is restored. Thus $\mathcal{L}_{2}=232454647871742$ becomes
$121343536768631$. Then apply the standard visualization, but move the distinguished yellow vertex to position $i$.\\
\indent Now, apply this convention to the actual sequences extracted from (\ref{sep-dev}). The result is shown in Fig. \ref{ratt}. What a relief: passing from sequence $1$ ($\mathcal{L}$) to sequence $2$ ($\mathcal{L}_{2}$), nothing changes except that the distinguished vertex moves by one place. In this light, the fact that the map $\mathscr{S}$ is well defined becomes tautological. Given a full triangulation of a regular $(n+1)$-gon, once all diagonals ending at one chosen vertex are removed, there is only one way to complete the triangulation again by adding diagonals starting from that very vertex.\\
\indent The same pattern continues up to sequence $8$: at each step the distinguished vertex simply advances by one. At sequence $9$ (which again begins at $1$), the whole triangulation rotates by one step in the positive direction, while the distinguished vertex jumps two vertices. If given some thought, it becomes clear why this is so.
\subsection*{Triangulation shifts}
What, then, is a glitch -- one of the shorter sequences -- in this language? How is it followed by a maximal sequence? Fig. \ref{ratt} can almost be read without words to see what is happening. At step $17$ (where the distinguished vertex is $8$), one edge ending at $8$ is deleted; at the same time, since the distinguished vertex changes from $8$ to $1$, the whole triangulation rotates one step forward. At step $18$ (distinguished vertex $1$), an edge ending at $1$ is added. Two steps later a similar event happens again: one edge disappears, and in the next step another is reinstated. A final observation makes the whole process easier to hold in the mind. When vertex $8$ is replaced by $1$, instead of imagining a rotation and a two-step jump, one may equally well imagine that there is no rotation at all, and that the distinguished vertex simply makes the usual one-step jump. After $24$ steps one then arrives at vertex $7$, from whose vantage point the very same initial triangulation is seen again.\\
\begin{figure*}
	\begin{center}
		\includegraphics[scale=0.45]{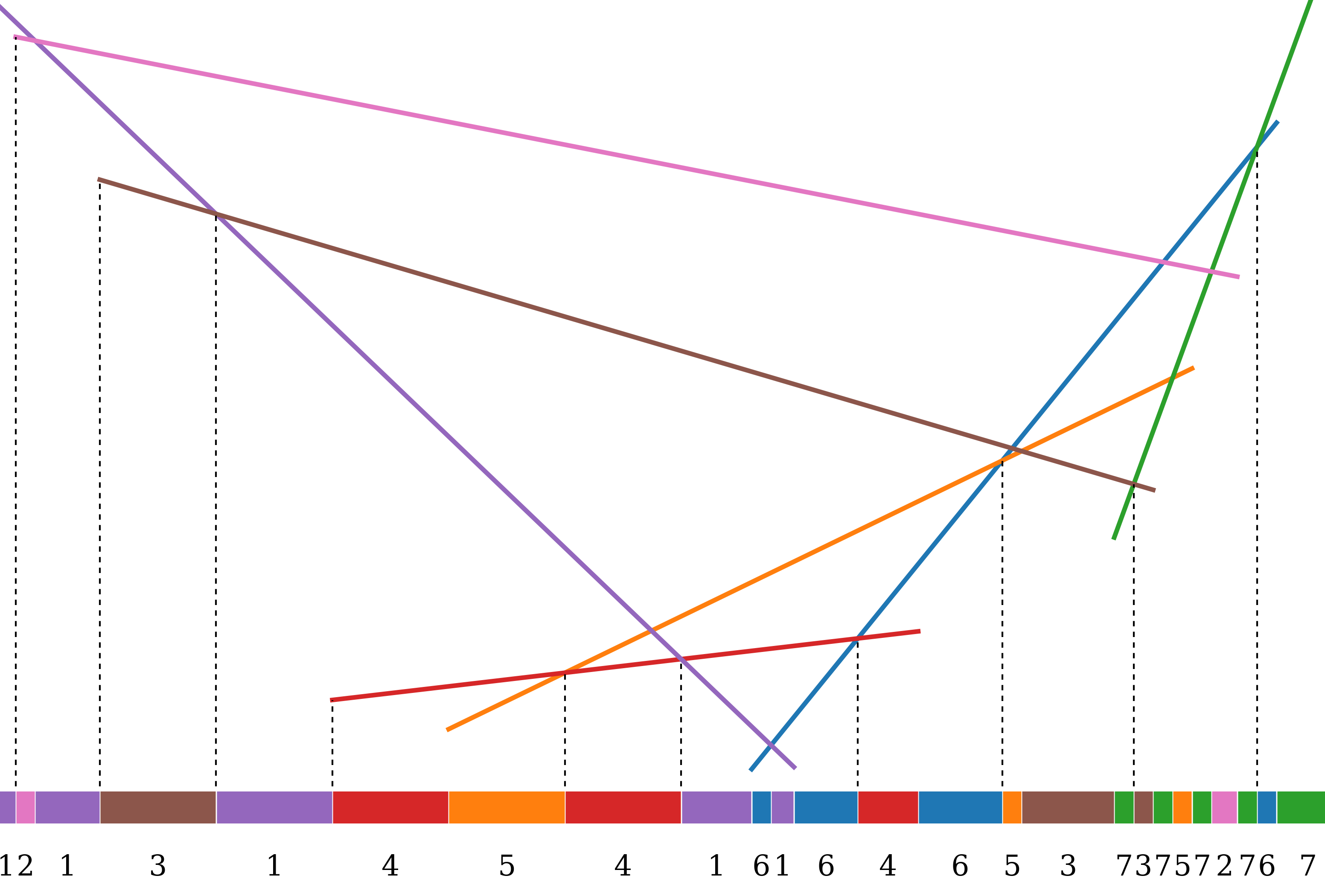} 
		\caption{The support for the lower envelope of $7$ segments with complexity $25$.}
		\label{support}
	\end{center}
\end{figure*}
\begin{figure*}
	\begin{center}
		\includegraphics[scale=0.86]{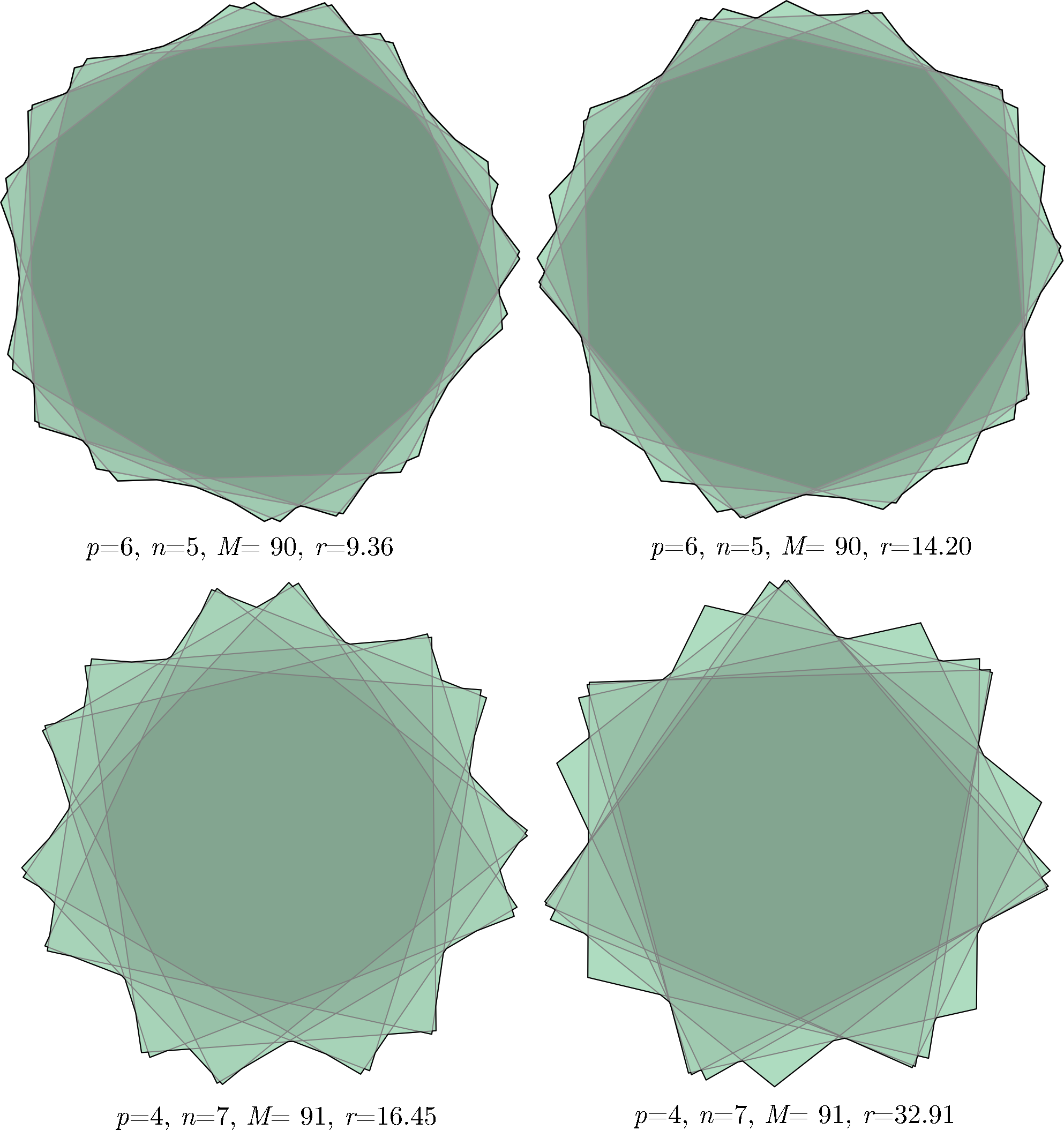} 
		\caption{Geometrically perfect polygon unions of inequivalent combinatorics. Top: $5$ hexagons (type 1). Bottom: $7$ quadrilaterals (type 2). The maximal ratio of edge-lengths is shown. This number is close to the optimal value.}
		\label{geom-perf}
	\end{center}
\end{figure*}
\indent The process can therefore be summarized in a compact form. Start from a distinguished vertex $1$ in a triangulation of a regular $(n+1)$-gon. Perform $3n$ moves. The $i$th move consists of removing a few diagonals at the distinguished vertex, inserting a few others there, counting their total number $t_i$, and then moving the distinguished vertex one step forward. The last vertex at which a move is performed is $n-2$ (after that one jumps to $n-1$ and stops). The goal is to return to the same triangulation, rotated by $-3$ vertices. The total combinatorial length of the corresponding triangle union is then $\nu=\sum_{i=1}^{3n}t_i+6n$.

\begin{figure}
	\begin{center}
		\includegraphics[scale=0.54]{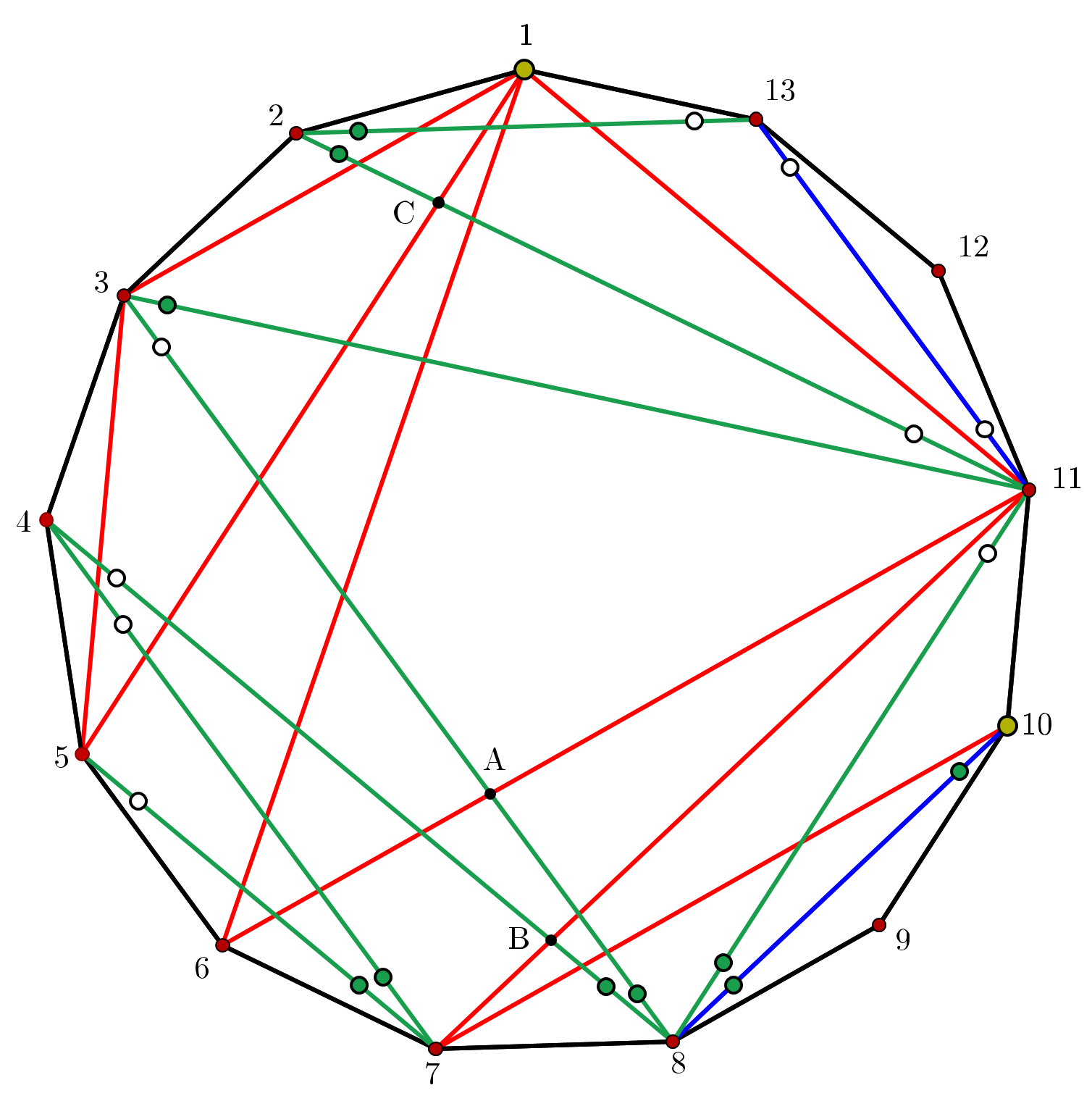} 
		\caption{The weight of this triangulation shift is $122=\mathscr{U}(12)-4$. The number of green dots (showing which tips of reinstated diagonals were counted in $10$ moves) is $10=\gamma(\mathcal{T}_{0})=\gamma(13)+4$. This triangulation gives a model for a triangle union which is $4$ edges short of the combinatorially maximal.}
		\label{aux}
	\end{center}
\end{figure}

\indent This is taken as the definition of a triangulation shift, and $\nu$ is called its weight. Intuitively, the largest weight is achieved by those initial triangulations $\mathcal{T}_0$ that overlap their own $-3$ rotation as much as possible. If no deletions were ever needed, then the triangulation would simply coincide with itself after rotation by $-3$. There are only two nontrivial examples of this phenomenon, namely when $n+1=6$ and $n+1=9$ (see Fig. \ref{gamma}); these exactly correspond to the Pentastar and the Octastar. In general, if no deletions were needed and the full $3$ turns are made, the weight would be $12n-12$. This suggests the auxiliary function $\gamma(\mathcal{T}_0)=12n-12-\nu(\mathcal{T}_0)$, measuring, so to speak, the defect of the triangulation shift. Now superimpose $\mathcal{T}_0$ with its own copy rotated by $-3$ turns. Colour the diagonals present only in the original in red, those present only in the rotated copy in green, and those common to both in blue. Fig. \ref{aux} shows an example. One starts from the yellow vertex $1$ and proceeds in the positive direction as far as the yellow vertex $10$. A red diagonal is deleted as soon as one of its endpoints is reached; a green diagonal is reinstated under the same rule, provided no red diagonal blocks it.
Any green diagonal of length $\geq 4$ (where the length of $(i,j)$ is defined as $\vert i-j\vert$) intersects its own red ``pre-image" -- as do diagonals $A$, $B$, and $C$ in Fig. \ref{aux}. Consequently, at least one endpoint of such a green diagonal will fail to be counted. Asymptotically, every triangulation contains at least $n/3$ such diagonals. To determine $\gamma$ exactly, brute-force computations were carried out in Python. Fig. \ref{gamma} shows the values for $5\leq n\leq 16$, together with one selected champion in each case (from $n+1=7$ onward, the champion is no longer unique). Once the pattern reveals itself, it becomes routine to prove that for $n\geq 5$,
$\gamma(n+1)=n+2-2\big{\lfloor}\frac{n+1}{3}\big{\rfloor}$.
\subsection*{Lower Envelopes}
Take $7$ linear univariate functions on the real line, labelled $1$ through $7$, and let $f$ be their pointwise minimum. If this minimum passes from one line to another, the abandoned line can never be jumped upon again. The \emph{support sequence} of $f$ -- the record of which function supplies the minimum as one moves along the $x$-axis -- therefore has length at most $7$ symbols.\\
\indent But if each linear function is defined only on its own interval, whose union is  connected, the support sequence may be longer. This is obvious. It is genuinely surprising that the support sequence can be considerably longer! For example, the $7$ segments in Fig. \ref{support} produce the support sequence
$1213145416164653737572767$, of length $25$. What is special about such a sequence? As with $\mathcal{L}$, no two consecutive symbols are equal. But here one immediately notices subsequences such as $1313$, $1414$, $3535$, $4646$, and so on. This is the crucial difference. It may be left to the reader to check that although a subsequence $abab$ may (and does) occur, none of the form $ababa$ can. That is exactly the defining property of a $(3,n)$-Davenport-Schinzel sequence.\\
\indent The opening part of the monograph by Sharir and Agarwal (\cite{sharir}, p. 12--42) is devoted to the study of $(3,n)$-DS sequences. The eventual answer -- that the maximal length is $2n\alpha(n)+O(n)$, due to Seth Pettie \cite{pettie} -- says little about the story of how that answer was discovered. In the book's preface one reads how Hart and Sharir first realized that superlinear growth is possible.\\
\indent There is another striking result there. Consider a maximal $(3,7)$-DS sequence of length $27$ (see OEIS A002004): $121314151616543272737475767$. It can be realized as the lower envelope of $7$ \emph{pseudosegments} (that is, connected pieces of $7$ pseudolines). Can it also be realized as the lower envelope of $7$ line segments? In other words, can the pseudosegment construction be stretched? That seems unknown. What is certain that there is little reason to believe -- and it is probably false; see \cite{sharir}, p. 112 -- that for all sufficiently large $n$, some maximal $(3,n)$-DS sequence can be realized as the lower envelope of $n$ segments. On the other hand, Wiernik and Sharir proved that there do exist $(3,n)$-DS sequences whose length differs asymptotically by at most a factor of $2$ from the true maximum and which do arise as lower envelopes of $n$ segments. As a mildly philosophical remark, this means that the Ackermann function, first encountered in logic and computability, appears in nature too -- provided that straight segments count as natural objects. Along with record-holding Turing machines \cite{strauss}, the aperiodic hat monotile, and tiles with large Heesch number, this is one of those phenomena in which mathematics transcends natural sciences and wanders far from everyday intuitions.\\
\indent To finish the construction, look again at Fig. \ref{support}. If each segment is completed to a triangle by giving it a left leg that leans slightly rightward and a right leg that leans slightly leftward, a simple polygon can be formed. Thus the problem of understanding $i(n)$ leads straight into deeper questions about unions and stretchability.
\begin{figure}
	\begin{center}
		\includegraphics[scale=0.48]{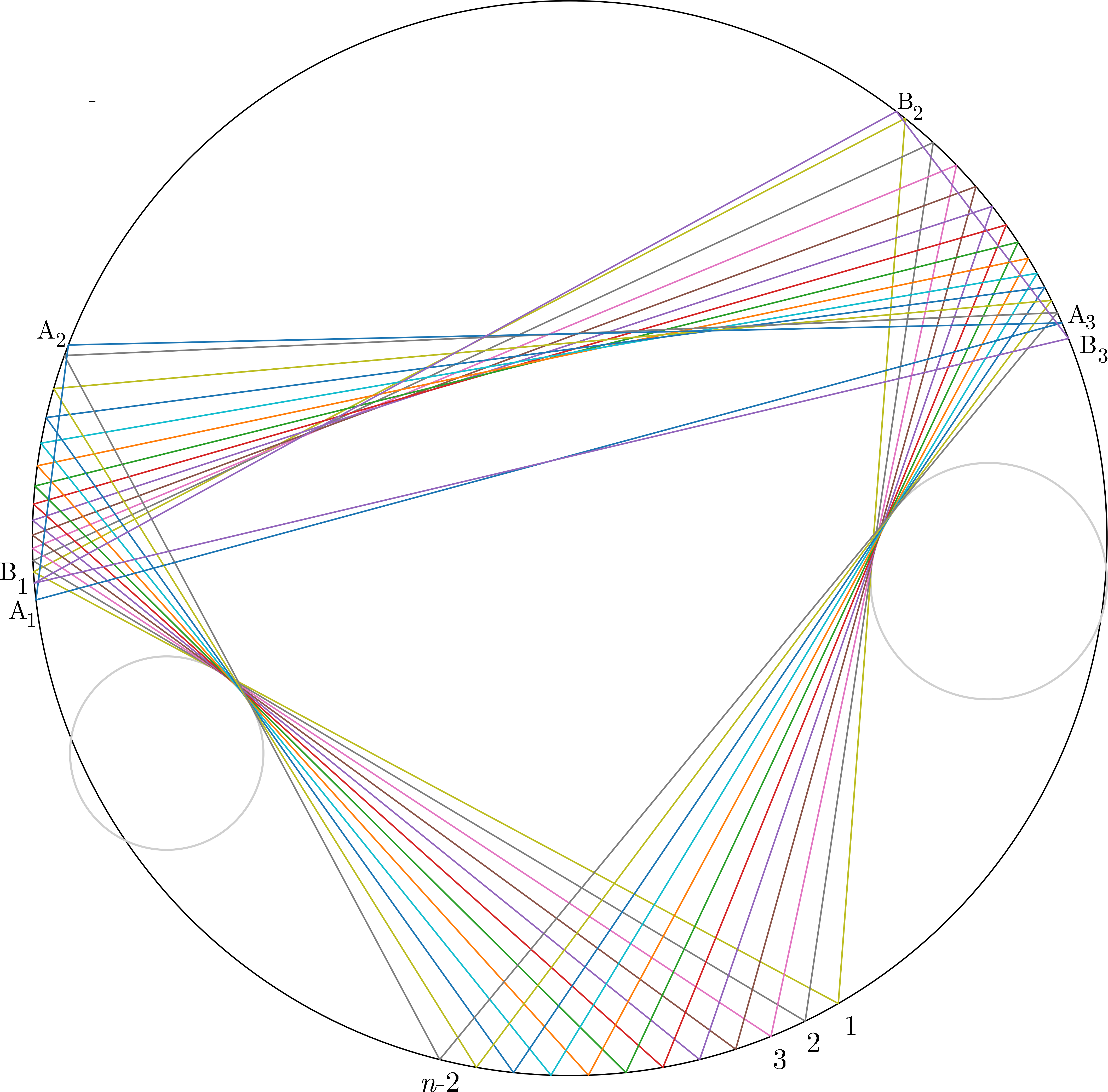} 
		\caption{A construction showing that $r(n)\geq 10n-10$. By moving $A_{2}$, $B_{2}$ down and $B_{3}$ up, a count $10n-7$ is achieved.}
		\label{two-cyc}
	\end{center}
\end{figure}
\subsection*{Open problems}
The sequences $c(n)$ (all vertices of all triangles on the unit circle) and $r(n)$ (with the additional requirement that the union be regular) were not treated in this paper. Since $r(n)\leq e(n)$, the latter is bounded above by a linear function. For $1\leq n\leq 8$ equality holds. This may be where a friendship ends: numerical search finds cyclically arranged $9$ triangles on the unit circle with $90$ edges, but none with $91$. On the other hand, Fig. \ref{two-cyc} yields a construction showing that $r(n)\geq 10n-7$. Some room for improvement still remains.\\
\indent At present, $c(n)$ is the most elusive of the four sequences. Of course
$r(n)\leq c(n)\leq i(n)$. Adding one triangle to Fig. \ref{two-cyc} can yield the lower bound $c(n)\geq 11n-12$. But beyond this, very little seems accessible. The lower-envelope method does not apply here. Indeed, Agarwal and Sharir explicitly pose the problem of determining the maximal complexity of the lower envelope of line segments whose two endpoints lie on a circular arc (\cite{sharir}, p. 112).\\
\indent Finally, triangulation shifts are by no means specific to triangle unions. They also help in studying $e_p(n)$, the maximal number of edges in a regular union of $n$ $p$-gons. One simply replaces $3n$ moves by $pn$ moves, and a rotation by $-3$ by a rotation by $-p$. The full problem is perhaps too broad, and is still premature. But it has one appealing subproblem of independent interest.\\
\indent A regular (and proper, meaning $n>1$) union of $n$ $p$-gons must satisfy
$(n+1)e_{p}(n)\leq (4n-2)pn$.
If equality holds, let such a union again be called perfect. Every perfect example must correspond to a triangulation of a regular $(n+1)$-gon that coincides with itself after rotation by $-p$ vertices. This leads to three basic types (for $p=3$, represented by the Star of David, the Pentastar, and the Octastar):
\begin{itemize}
	\item[]{Type $1$}: $n=p-1$, $M=2(p-1)(2p-3)$.
	\item[]{Type $2$}: $n=2p-1$, $M=(2p-1)(4p-3)$. 
	\item[]{Type $3$}: $n=3p-1$, $M=2(2p-1)(3p-1)$.
\end{itemize}
If such a perfect union can be realized with straight lines and also possesses an $n$-fold rotational symmetry, let it be called \emph{geometrically perfect}. For $p=6$, $n=5$ (Type $1$), and for $p=4$, $n=7$ (Type $2$), there are two combinatorially distinct perfect unions. This reflects the fact that there are two inequivalent triangulations of a regular hexagon (up to the action of the dihedral group $\mathbb{D}_6$), and likewise two inequivalent center-symmetric triangulations of a regular octagon (up to $\mathbb{D}_8$). All four can be realized geometrically (Fig. \ref{geom-perf}). In the last example, however, the shortest edge must be very short indeed. It is plausible that only finitely many geometrically perfect proper polygon unions exist. If so, they would make a slightly eccentric addition to the other finite catalogues of geometry ($5$ Platonic solids, $11$ Archimedean tilings, $15$ families of pentagonal tilings, and so on).\\
\indent To conclude, the torch is gladly passed to any reader inclined to take it up. The first four questions below look especially approachable.
\begin{itemize}
	\item[$1$)]Find $n$ for which $i(n)>e(n)$.
	\item[$2$)]Prove rigorously that $r(9)=90$.
	\item[$3$)]Improve the lower bound for $r(n)$.
	\item[$4$)]Find a perfect union which is not realisable as geometrically perfect.
	\item[$5$)]What arguments support, and what arguments oppose, the possible linear growth of $c(n)$?
	\item[$6$)]Can $e(n+1)-e(n)=14$ occur?
	\item[$7$)] Find the value of the constant 
	\begin{eqnarray*}
		\limsup\limits_{n\rightarrow\infty}\frac{e(n)}{n}\in\bigg{[}11,\frac{35}{3}\bigg{]}.
	\end{eqnarray*}
\end{itemize}
Since $\mathscr{U}(n+1)-\mathscr{U}(n)$, from $n=5$ onward, repeats periodically as $11,11,13$, the penultimate question may well have an affirmative answer. As for the last one, intuition remains, sadly, silent. 

\subsection*{Funding declaration}
The paper received no external funding.


\begin{thebibliography}{2}
\bibitem{oeis}{\it The On-Line Encyclopedia of Integer Sequences}, sequences A002004, A006066, A090338,  A090339, A241600, A250001, A375986.

\bibitem{ackerman} {\sc E. Ackerman, B. Keszegh, Rote G.}, An almost optimal bound on the number of intersections of two simple polygons. {\it Discrete Comput. Geom.} {\bf 68}(4) (2022), 1049--1077, \url{https://link.springer.com/article/10.1007/s00454-022-00438-0}.


\bibitem{agarwal}{\sc P.K.  Agarwal, J. Pach, M. Sharir}, State of the union (of geometric objects). {\it Surveys on discrete and computational geometry}, 9--48, Contemp. Math., 453, Amer. Math. Soc., 2008.

\bibitem{horizon}{\sc M. Bern, D. Eppstein, P. Plassmann, F. Yao}, Horizon theorems for lines and polygons. {\it Discrete and computational geometry}, 45–66, DIMACS Ser. Discrete Math. Theoret. Comput. Sci., 6, Amer. Math. Soc,, RI (1991).	
	
\bibitem{zone}{\sc B. Chazelle, L.J. Guibas, D.T. Lee}, The power of geometric duality. {\it BIT} {\bf 25}(1) (1985), 76-–90.	

\bibitem{cutler}{\sc D.O.H. Cutler, N.J.A. Sloane}, Cutting a Pancake with an Exotic Knife, \url{https://arxiv.org/pdf/2511.15864}.

\bibitem{dallant}{\sc J. Dallant}, Improved Bound on the Number of Pseudoline Arrangements via the Zone Theorem, 2025, \url{https://arxiv.org/abs/2502.20909}.

\bibitem{edelsbrunner}{\sc H. Edelsbrunner, L. Guibas, J. Pach, R. Pollack, R. Seidel, M. Sharir}, Arrangements of curves in the plane—topology, combinatorics, and algorithms. {\it Theoret. Comput. Sci.} {\bf 92}(2) (1992), 319--336.

\bibitem{felsner}{\sc S. Felsner, J.E. Goodman}, Pseudo-line arrangements, {\it Handbook of Discrete and Computational Geometry}, Taylor and Francis, 2017. \url{https://doi.org/10.1201/9781315119601}

\bibitem{strauss}{\sc Ch. Goodman-Strauss}, Can't Decide? Undecide! {\it Notices Amer. Math. Soc.}, {\bf 57} (3) (2010), 343–-356, \url{https://www.ams.org/notices/201003/rtx100300343p.pdf}

\bibitem{peled}{\sc S. Har-Peled}, {\it Union Complexity of Pseudodiscs}, March 21, 2023, \url{https://courses.grainger.illinois.edu/cs498sh3/sp2023/lec/18/union_pseudodisks.pdf}.
	
\bibitem{fat}{\sc J. Matou\v{s}ek, J. Pach, M. Sharir, Sh. Sifrony, E. Welzl}, Fat triangles determine linearly many holes. {\it SIAM J. Comput.} {\bf 23}(1) (1994), 154--169, \url{https://epubs.siam.org/doi/10.1137/S009753979018330X}

\bibitem{pettie}{\sc S. Pettie}, Sharp bounds on Davenport-Schinzel sequences of every order. {\it J. ACM} {\bf 62}(5) (2015), Art. 36, 40 pp.

\bibitem{pinchasi}{\sc R. Pinchasi}, The zone theorem revisited. Manuscript, 2011. \url{https://sites.google.com/view/homepage-of-rom-pinchasi/home}.
	
\bibitem{roselle}{\sc D.P. Roselle}, An algorithmic approach to Davenport-Schinzel sequences. {\it Utilitas Math.} 6 (1974), 91--93.
	
\bibitem{schwartz}{\sc R.E. Schwartz}, Two geometric pictures of Farey addition. {\it Math. Intelligencer} {\bf 47}(3) (2025), no. 3, 271--275. \url{https://link.springer.com/article/10.1007/s00283-025-10410-4}

\bibitem{sharir}{\sc M. Sharir, P.K. Agarwal}, Davenport-Schinzel sequences and their geometric applications. {\it Cambridge University Press}, Cambridge, 1995.
	
\bibitem{sloane}{\sc N.J.A. Sloane}, The On-Line Encyclopedia of Integer Sequences. {\it Notices Amer. Math. Soc.} {\bf 65}(9),  1062--1074 (2018).

\bibitem{suvorov}{\sc P. Suvorov}, Isotopic but not rigidly isotopic plane systems of straight lines. In: Topology and Geometry — Rohlin Seminar. Lecture Notes in Mathematics, vol 1346. (1988) Springer. \url{https://doi.org/10.1007/BFb0082793}.
\end{thebibliography}
\end{document}